\documentclass[11pt]{article}
\usepackage{amssymb}
\usepackage{amsmath}
\usepackage{theorem}
\usepackage{algorithm}
\usepackage{graphicx}
\usepackage[usenames]{color}
\usepackage{multirow}

\newcommand{\R}{\mathbb{R}}

\newcommand{\ve}[1]{ #1}
\newtheorem{Thm}{Theorem}[section]
\newtheorem{Cor}{Corollary}[section]
\newtheorem{Lemma}{Lemma}[section]

  {\end{list}}

\theoremheaderfont{\scshape}
\makeatletter

  \gdef\listctr{list\romannumeral\the\@listdepth}\expandafter
  
\makeatother

\newenvironment{AlgorithmSteps}[1][1]{%
  \begin{list}{\csname label\listctr\endcsname}{%
      \usecounter{\listctr}
      
      \settowidth{\labelwidth}{\textsc{Step\ #1.}}%
      \setlength{\leftmargin}{\labelwidth}\addtolength{\leftmargin}{\labelsep}}}%
  {\end{list}}

\def\x{ {\ve x}}
\def\y{ {\ve y}}
\def\u{ {\ve u}}
\def\w{{\ve w}}
\def\z{ {\ve z}}

\def\xk{ \x^{(k)}}
\def\yk{ \y^{(k)}}
\def\tildeyk{ \tilde{\y}^{(k)}}
\def\uk{ \u^{(k)}}
\def\zk{ \z^{(k)}}
\def\yk{ \y^{(k)}}
\def\gk{ \w^{(k)}}
\def\xkk{ \x^{(k+1)}}

\def\ak{{\alpha_k}}
\def\ek{{\epsilon_k}}
\def\lmin{\lambda_{min}}\def\lmax{\lambda_{max}}
\def\ykk{ \y^{(k+1)}}

\def\xtilde{\tilde\x}\def\kinN{{}}
\def\diam{\mbox{diam}}
\def\dom{\mbox{dom}}
\def\argmin{\mbox{argmin}}
\def\R{\mathbb R}
\def\frec{f^{rec}}
\def\flev{f^{lev}}

\def\endproof{\hfill$\square$\vspace{0.3cm}\\}
\newcommand\silviacorr{}
\newcommand\valeriacorr{}
\newcommand\proof{{\it Proof. }}

\begin{document}
\title{Scaling techniques for $\epsilon$--subgradient projection methods\thanks{This work has been partially supported by MIUR (Italian Ministry for University and Research), under the projects FIRB - Futuro in Ricerca 2012, contract RBFR12M3AC, and by the Italian Spinner2013 PhD project {\it High-complexity inverse problems in biomedical applications and social systems}. The Italian GNCS - INdAM (Gruppo Nazionale per il Calcolo Scientifico - Istituto Nazionale di Alta Matematica) is also acknowledged.}}
\author{S. Bonettini\footnotemark[2], A. Benfenati\footnotemark[3], V. Ruggiero\footnotemark[4]}
\maketitle
\renewcommand{\thefootnote}{\fnsymbol{footnote}}

\footnotetext[2]{Dipartimento di Matematica e di Informatica, Universit\`a di Ferrara, Polo Scientifico Tecnologico, Blocco B, Via Saragat 1, I-44122 Ferrara, Italy ({\tt silvia.bonettini@unife.it}).}
\footnotetext[3]{Dipartimento di Matematica e di Informatica, Universit\`a di Ferrara, Polo Scientifico Tecnologico, Blocco B, Via Saragat 1, I-44122 Ferrara, Italy ({\tt alessandro.benfenati@unife.it}).}
\footnotetext[4]{Dipartimento di Matematica e di Informatica, Universit\`a di Ferrara, Polo Scientifico Tecnologico, Blocco B, Via Saragat 1, I-44122 Ferrara, Italy ({\tt valeria.ruggiero@unife.it}).}
\renewcommand{\thefootnote}{\arabic{footnote}}

\begin{abstract}\noindent
The recent literature on first order methods for smooth optimization shows that significant improvements on the practical convergence behaviour can be achieved with variable stepsize and scaling for the gradient, making this class of algorithms attractive for a variety of relevant applications. In this paper we introduce a variable metric in the context of the $\epsilon$-subgradient projection methods for nonsmooth, constrained, convex problems, in combination with two different stepsize selection strategies. We develop the theoretical convergence analysis of the proposed approach and we also discuss practical implementation issues, as the choice of the scaling matrix. In order to illustrate the effectiveness of the method, we consider a specific problem in the image restoration framework and we numerically evaluate the effects of a variable scaling and of the steplength selection strategy on the convergence behaviour.
\end{abstract}
%

\pagestyle{myheadings}
\thispagestyle{plain}
\markboth{S. BONETTINI, A. BENFENATI, V. RUGGIERO}{A SCALED $\epsilon$-SUBGRADIENT PROJECTION METHOD}
\section{Introduction}
In this paper we consider the problem
\begin{equation}
\min_{x\in X} f(x)
\label{minf}
\end{equation}
where \silviacorr{$f: \R^n\rightarrow \R\cup \{\infty\}$ } is a convex, proper, lower semicontinuous function and $X$ is a nonempty, closed, convex subset of $\R^n$ \silviacorr{contained in the domain of $f$}. We denote by $X^*$ the set of the solutions of \eqref{minf}. We are interested in the case where $f$ is nondifferentiable and a subgradient or an approximate subgradient of $f$ can be easily computed.
This arises for example in duality and minimax contexts. A well-known method to solve the problem (\ref{minf}) is the $\epsilon$-subgradient method:
\begin{equation}
\label{subgradient}
\xkk = P_{X}\left(\xk-\ak \uk\right)
\end{equation}
where $\uk\in\partial_{\ek} f(\xk)$ for some $\ek\geq 0$,  $\alpha_k$ is a positive stepsize and $P_X(\cdot)$ is the Euclidean projection operator on the set $X$. The choice $\ek=0$ for all $k$ corresponds to the subgradient method, which has been extensively investigated (see for example the contributions collected in \cite{Ermoliev83,Goffin77,Polyak1987,Shor85}).

The more general case allowing $\ek>0$ was introduced and developed in \cite{Ermoliev83,Polyak1987}, while more recent convergence results under different assumptions are given in \cite{Alberetal98,Correa93,Kiwiel-2004,Larsson03,Robinson99,Solodov1998}. A typical assumption on the sequence $\{\epsilon_k\}$ is that
\begin{equation}\label{epsilontozero}\lim_{k\rightarrow \infty}\ek = 0 \end{equation}
and, in this case, the subgradient and the $\epsilon$-subgradient methods have very similar convergence properties. In the following discussion we assume that \eqref{epsilontozero} holds.

The $\epsilon$-subgradient method is interesting in itself, since, when the projection on $X$ is easy to compute and an approximate subgradient is available, it can be easily implemented, but it is also a useful tool to analyze the theoretical convergence properties of a variety of algorithms \cite{BonettiniRuggiero12,Esseretal09,Neto09,Robinson99}.

We can distinguish different variants of the method \eqref{subgradient} according to the rule adopted to select the stepsize. We list below the most studied stepsize choices for subgradient methods, which, with minor modifications, could be applied also to the case $\epsilon_k>0$:
\begin{itemize}
\item[$({\mathcal R}_1)$] the \emph{constant stepsize} rule $\ak = \alpha > 0$;
\item[$({\mathcal R}_2)$] the \emph{Polyak rule}
\begin{equation}\nonumber
\ak = c_k\frac{f(\xk)-f^*}{\|\uk\|^2}\ \  \mbox{ or }\ \  \ak = c_k\frac{f(\xk)-f^*}{\max\{1,\|\uk\|^2\}}\ \ \ c_k\in(0,2)
\end{equation}
where $f^*$ is the minimum of $f$
\item[$({\mathcal R}_3)$] the \emph{Ermoliev} or \emph{diminishing, divergent series} stepsize rule, which includes any sequences $\{\ak\}_\kinN$ such that
\begin{equation}\label{Ermoliev}
\ak > 0\ \ \lim_{k\rightarrow\infty}\ak = 0 \ \ \sum_{k=0}^{\infty} \ak = \infty
\end{equation}
\item[$({\mathcal R}_4)$] the \emph{diminishing, divergent series, square summable} stepsize rule, which, in addition to \eqref{Ermoliev}, also requires $\sum_{k=0}^{\infty} \ak^2< \infty$.
\item[$({\mathcal R}_5)$] the \emph{dynamic} or \emph{adaptive} stepsize rule
\begin{equation}\label{dynamic}
\ak = \frac{f(\xk)-f_k}{\|\uk\|^2}\ \  \mbox{ or }\ \ \ak = \frac{f(\xk)-f_k}{\max\{1,\|\uk\|^2\}}
\end{equation}
where $f_k$ is an adaptively computed estimate of $f^*$; several further variants of this rule, which can be considered as an approximation of $({\mathcal R}_2)$ when $f^*$ is not known, depend on how $f_k$ is defined.
\end{itemize}
Keeping a constant stepsize as in $({\mathcal R}_1)$, only the convergence of a subsequence of $\{f(\xk)\}_\kinN$ to a possibly suboptimal value is established, i.e. $\lim\inf_k f(\xk)\leq f^* + C\alpha$, for some positive constant $C$ \cite{Bertsekas2012,Nedic01}; for rules $({\mathcal R}_3)$ and $({\mathcal R}_5)$ stronger results have been proved \cite{Goffin77,Larsson03,Nedic01}, showing that $\lim_k f(\xk) = f^*$ and $\min_{\x^*\in X^*}\|\x^*-\xk\|\rightarrow 0$ (with the assumption $\epsilon_k = 0$ for the latter case). Finally, the convergence of the sequence $\{\xk\}_\kinN$ to a solution of \eqref{minf} can be proved in the cases $({\mathcal R}_2)$, with $\ek = 0$, and $({\mathcal R}_4)$ \cite{Alberetal98,Larsson03}.

Results about an optimal stepsize choice to obtain a suboptimal rate of convergence for the subgradient method are reported in \cite{Nesterov04}; in \cite[Sec. 6.3]{Bertsekas2012}, a convergence analysis is performed for different stepsize choices while in \cite{Auslander09} analogous results are obtained with respect to non Euclidean metrics.

The key property that the stepsize parameter has to induce on the iterates \eqref{subgradient} which is exploited in the standard convergence analysis for subgradient methods is the \emph{quasi-F\'{e}jer monotonicity} with respect to the set $X^*$
$$\|\xkk-\x^*\|^2\leq \|\xk-\x^*\|^2 + \eta_k \ \ \ \forall x^*\in X^*$$
for some nonnegative sequence $\{\eta_k\}_\kinN$ such that $\sum \eta_k<\infty$ (see \cite{Baushke-Combettes-2011,Combettes-2001}).

It is worth noticing that the stepsize in subgradient methods plays quite a different role than in the smooth case, where an analogous parameter is employed to ensure the sufficient decrease of the objective function, and, in some kind of schemes, also to accelerate the convergence, for example by means of the well known Barzilai--Borwein rules \cite{Barzilai1988,Dai2006a,Frassoldati2008} (see also \cite{DeAsmundis-etal-2013, Fletcher2012} for recent developments in this field).
Thus, these valid approaches to the stepsize selection for the smooth case are difficult to extend to the method \eqref{subgradient}.

On the other side, recent advances in the context of gradient based methods show that introducing a variable scaling matrix for the gradient can lead to significant improvements on the practical performances \cite{Bardsley2006,Bonettini2009,Bonettini2013,Hager2009,Theys-etal-2009}, especially on large scale and ill conditioned problems. \silviacorr{Variable metric was introduced also in \cite{Combettes-Vu-2014} in the context of monotone operators and in \cite[Chapter 5]{Nedic-2002} in the subgradient methods for unconstrained optimization.}

Motivated by this,
we propose to introduce a variable scaling matrix for the $\epsilon$-subgradient vector in the iteration \eqref{subgradient}.
More precisely, the contribution of this paper is to provide the convergence analysis, under standard assumptions, of the following scaled $\epsilon$-subgradient scheme:
\begin{equation}\label{method}
\xkk = P_{X,D_k^{-1}}\left(\xk-\ak D_k{\uk}\right)
\end{equation}
where $D_k$ is a symmetric positive definite matrix with bounded eigenvalues, the projection operator is defined as
\begin{equation}\label{proj}
P_{X,D_k^{-1}}(\x) = \mbox{argmin}_{\z\in X} 
\ (\x-\z)D_k^{-1}(\x-\z)
\end{equation}
and $\ak$ is chosen either as an {\it a priori} selected sequence obeying to the diminishing, divergent series, summable squares stepsize rule $({\mathcal R}_4)$, or with an adaptive rule $({\mathcal R}_5)$ of Br\"{a}nnlund's type \cite{Brannlund95,Goffin99,Nedic01}.

In the first case, assuming that the set $X^*$ is nonempty, we prove the convergence of the sequence $\{\xk\}_{k\in \mathbb N}$ to a point $x^*\in X^*$, while in latter one we prove the convergence of the sequence $\{f(\xk)\}_{k\in \mathbb N}$ to the minimum value $f^*$.

A further contribution of the paper is to introduce, as special case of \eqref{method} and also as a generalization of the method in \cite{BonettiniRuggiero12}, a Scaled Primal--Dual Hybrid Gradient (SPDHG) method, which applies to the case
\begin{equation}\label{minf0f1}
\min_{x\in X} f_0(x)+f_1(Ax)
\end{equation}
where $f_0$ and $f_1$ are convex, proper, lower semicontinuous functions and $A$ is a linear operator. Many relevant problems can be modeled, as, for example, the restoration of images approached from the Bayesian paradigm \cite{Geman-Geman-1984}.

\silviacorr{Several methods for the solution of \eqref{minf0f1} have been developed in the recent literature: in particular, when $f_0$ is continuously differentiable with Lipschitz continuous gradient, suitable splitting and forward-backward methods can be applied \cite{Combettes-Vu-2014,Condat-2013,Lorenz-Pock-2014}. As we will show in Section \ref{sec4}, problem \eqref{minf0f1} can be handled also by SPDHG, even when $f_0$ is non differentiable or its gradient is not Lipschitz continuous on $\dom (f_0)$.}
In particular, we provide an especially tailored implementation of SPDHG for the
Total Variation (TV) restoration of images corrupted by Poisson noise. This problem is related to several applications such as astronomical imaging, electronic microscopy, single particle emission computed tomography
(SPECT) and positron emission tomography (PET), and a variety of specialized methods have been proposed for its solution (see \cite{Bardsley09,BonettiniRuggiero2011,Dupé-etal-2011a,Dupé-etal-2011b,Figueiredo2010,Setzer10} and references therein). For this special case of SPDHG, we devise an effective strategy to choose the scaling matrix, showing that significant improvements on the practical convergence speed can be obtained.

The paper is organized as follows. In Section \ref{sec2} we present a convergence analysis for the scaled $\epsilon$-subgradient method \eqref{method} when the stepsizes are chosen according to a diminishing, divergent series, square summable stepsize rule. At the end of the section our results are also compared to the very recent works \cite{Combettes-Vu-2013,Combettes-Vu-2014}, where variable metrics are studied from the point of view of more general operators.

Building on this material, in Section \ref{sec3} we propose a generalization to variable scaling and approximate subgradients of the \emph{level algorithm} in \cite{Goffin99,Nedic01}, based on a dynamic stepsize selection rule, showing that in our more general settings the main properties still hold. Further, in Section \ref{sec4} we consider problem \eqref{minf0f1} and we present the SPDHG method, proving its convergence as special case of an $\epsilon$-subgradient scheme. In order to illustrate a practical implementation of SPDHG, in Section \ref{sec5} we describe the problem of deblurring an image corrupted by Poisson noise via the TV regularization. A suitable scaling for the SPDHG method is discussed and an algorithm for its computation is detailed.
In Section \ref{sec6}, we describe some numerical simulations concerning the considered application, with the aim to evaluate the effectiveness of the scaling technique in the $\epsilon$-subgradient methods in combination with the two stepsize selection strategies analyzed in the previous sections. The numerical experiments show that a suitable selection of the scaling matrix can be a very effective tool to improve the convergence behaviour also in nonsmooth methods.
Finally, some concluding remarks are given in Section \ref{sec7}.
\paragraph{Notations and definitions.}
In the following, $\|\cdot\|$ denotes the Euclidean vector or matrix norm. Given $\x\in \R^n$ and a symmetric and positive definite matrix $D$ of order $n$, $\|\x\|_D$ denotes the energy norm, i.e. $\|\x\|_D=\sqrt{\x^T D \x}$, and $P_{X,D}(x)=\argmin_{y\in X}\|y-x\|_{D}^2$. By $\dom (f)$ we indicate the domain of any function \silviacorr{$f:\R^n\rightarrow \R\cup\{\infty\}$}, i.e. $\dom(f) = \{\x\in \R^n: f(\x)<\infty\}$, while $\diam(X)$ denotes the diameter of the closed, convex set $X\subset \R^n$, $\diam(X) = \max_{\x,\z\in X}\|\x-\z\|$. \silviacorr{The Fenchel dual or conjugate of $f$ is defined as $f^*(\y) = \sup_{\x\in \R^n} \x^T\y-f(\x)$.} Furthermore, we recall that the $\epsilon$-subdifferential of $f$ at $x\in \dom (f)$ for some $\epsilon\in \R$, $\epsilon \geq 0$, is the set
$$\partial_{\epsilon} f(x)= \{p\in \R^n:\ f(z)\geq f(x) + p^T(z-x) -\epsilon,\ \ \forall z\in \R^n\} $$
Any vector $p\in \partial_{\epsilon} f(x)$ is an $\epsilon$-subgradient of $f$ at $x$ \cite[\S 23]{Rokafellar70}.\\
\silviacorr{If $f(x) = \sum_{i=1}^n \beta_if_i(x)$, where $\beta_i\geq 0$, $u_i\in \partial_{\epsilon_i}f_i(x)$ and $x\in \bigcap_{i=1}^n \dom(f_i)$, then $\sum_{i=1}^n \beta_iu_i\in \partial_{\epsilon}f(x)$, where $\epsilon = \sum_{i=1}^n \epsilon_i$.
The proof of this property can be found in \cite{BonettiniRuggiero12,Rokafellar70}.}
\section{Convergence analysis with square summable stepsize sequences}\label{sec2}

In this section we show that with the rule $({\mathcal R}_4)$ the method \eqref{method}  generates a sequence of points converging to a solution of \eqref{minf}, under standard assumptions on
the error sequence $\ek$ and the scaling matrices $D_k$.

Before giving the main convergence result, we prove the following two useful lemmata concerning the energy norm, the corresponding projection operator and the summable sequences.
\begin{Lemma}\label{lemma:0}
Let $D$ be a symmetric positive definite matrix.
The following relations hold.
\begin{itemize}
\item[(i)] For all $\ve a, b, c\in \R^n$, we have
\begin{equation*}
\|a - b\|_D^2 + \|b - c\|_D^2 - \|a-c\|^2_D = 2 (b-a)^TD(b-c)
\end{equation*}
\item [(ii)] For any $\x\in \R^n$, $\z\in X$ we have
\begin{eqnarray}
(P_{X,D^{-1}}(\x)-\x)^T D^{-1} (\z-P_{X,D^{-1}}(\x))&\geq &0 \label{dis1}\\
(P_{X,D^{-1}}(\x)-\x)^T D^{-1} (\z-\x)&\geq &0 \label{dis2}
\end{eqnarray}
\item[(iii)] For any $\x, \z \in \R^n$ we have
\begin{equation*}
    \|P_{X,D^{-1}}(\x) - P_{X,D^{-1}}(\z)\|_{D^{-1}}\leq \|\x-\z\|_{D^{-1}}
    \end{equation*}
\item[(iv)] Let $L$ a positive number such that $\|D\|=\lambda_{max}(D)\leq L$ and $\|D^{-1}\|=\frac{1}{\lambda_{min}(D)}\leq L$, then
    \valeriacorr{
    \begin{equation}\label{2.2a}
    \|P_{X,D^{-1}}(\x) - P_{X,D^{-1}}(\z)\|\leq L \|\x-\z\|
    \end{equation}}
\end{itemize}
\end{Lemma}
{\it Proof.} Part (i) directly follows from the definition of $\|\cdot\|_D$, while the optimality conditions of the minimum problem \eqref{proj} yield part (ii) (see also \cite[Proposition 3.7]{Bertsekas-Tsitsiklis-1989}). Part (iii) is a consequence of \eqref{dis1} and of the Cauchy--Schwartz inequality. For (iv), see \cite[Lemma 2.1]{Bonettini2009}.
\endproof
\begin{Lemma}\label{lemma:3}
Let $\{L_k\}$ be a sequence of positive numbers such that $L_k^2 = 1+\gamma_k$, $\gamma_k \geq 0$, where $\sum_{k=0}^\infty \gamma_k<\infty$. Let $\theta_k = \prod_{j=0}^kL_j^2$ for any $k\geq 0$. Then the sequence $\{\theta_k \}$ is  bounded.
\end{Lemma}
{\it Proof.} We want to show that there exists a constant $M>0$ such that $\theta_k\leq M$ for all $k\geq 0$. By the monotonicity of the logarithm, this is true if and only if $\log(\theta_k)\leq \log(M)$. By definition of $\theta_k$ we have
\begin{equation}\label{serie}
\log(\theta_k) = \sum_{j=0}^k \log(L_j^2)\leq\sum_{j=0}^\infty \log(L_j^2)
\end{equation}
Thus, if the series on the right hand side of \eqref{serie} converges, the quantities $\theta_k$ are bounded for all $k$. We observe that, since $L_k^2 = 1+\gamma_k$, where $\gamma_k\rightarrow 0$ as $k$ diverges,  by the known limit $\lim_{k\rightarrow\infty} \frac{\log(1+\gamma_k)}{\gamma_k} = 1$, the series $\sum_{j=0}^\infty \log(L_j^2)$ and $\sum_{j=0}^\infty \gamma_j$ have the same behavior. Thus, since  by hypothesis the latter one is convergent, the theorem follows.\endproof
We are now ready to present the main convergence result about the method \eqref{method} with diminishing, divergent series, summable square stepsize sequences, whose proof is developed using similar techniques as in \cite[Lemma 1]{Alberetal98} and \cite[Theorem 8]{Larsson03}.
\begin{Thm}\label{teo:1}
Let $\{\xk\}\subset \Omega$ be the sequence generated by iteration \eqref{method}, where $\uk\in \partial_\ek f(\xk)$, for a given sequence $\{\ek\}\subset \R$, $\ek \geq 0$.
Assume that the set of the solutions of \eqref{minf} $X^*$ is nonempty and that there exists a positive constant $\rho$ such that $\|\uk\|\leq \rho$ and a sequence of positive numbers $\{L_k\}$ such that $\|D_k\|\leq L_k$, $\|D_k^{-1}\|\leq L_k$, with $1\leq L_k\leq L$ for some positive constant $L$, for all $k\geq 0$. If the following conditions holds
\begin{eqnarray}
&\displaystyle
\lim_{k\to\infty}\epsilon_k = 0\label{hp_ezerovera}   \\
&\displaystyle\sum_{k=0}^\infty\ak = \infty\label{hp_aezero} \\
&\ \ \ \displaystyle\sum_{k=0}^\infty \ak^2<\infty\ \ \ \sum_{k=0}^\infty \epsilon_k\ak<\infty\label{hp_sumepsilon}\\
&\displaystyle L_k^2 = 1 + \gamma_k\ \ \ \sum_{k=0}^{\infty} \gamma_k < \infty\label{hp_gammak}
\end{eqnarray}
then, the sequence $\{\xk\}$ converges to a solution of \eqref{minf}.
\end{Thm}
{\it Proof.} For all $k$ let us define $\zk = \xk-\ak D_k\uk$. By Lemma \ref{lemma:0} \valeriacorr{part (iii) we have that
\begin{equation}\label{ine0}
\|\xkk-\xk\|_{D_k^{-1}}\leq \|\zk-\xk\|_{D_k^{-1}}\leq L_k^{\frac 1 2}\ak \|\uk\| \leq L_k^{\frac 1 2}\ak \rho
\end{equation}}
Then, thanks to Lemma \ref{lemma:0}, part (i), for any $\xtilde\in X^*$ we can write
\begin{eqnarray}
&\valeriacorr{L_k}\ak^2\rho^2 &+\|\xk-\xtilde\|_{D_k^{-1}}^2-\|\xkk-\xtilde\|_{D_k^{-1}}^2\geq\nonumber\\
&&\geq\|\xkk-\xk\|_{D_k^{-1}}^2+\|\xk-\xtilde\|_{D_k^{-1}}^2-\|\xkk-\xtilde\|_{D_k^{-1}}^2\nonumber\\\nonumber
&&= 2( \xk-\xtilde)^TD_k^{-1}(\xk-\xkk) \\\nonumber
&& = 2(\xk-\xtilde)^TD_k^{-1}(\xk-\zk) + 2 (\xk-\xtilde)^TD_k^{-1}(\zk-\xkk)\\\nonumber
&& = 2\ak (\xk-\xtilde)^T \uk+ 2(\xk-\zk)^TD_k^{-1}(\zk-\xkk) + 2 (\zk-\xtilde)^TD_k^{-1}(\zk-\xkk)\\\nonumber
&& = 2\ak (\xk-\xtilde)^T\uk+ 2(\xk-\zk)^TD_k^{-1}(\zk-\xkk) +\nonumber\\ &&\ \ \ + 2 (\zk-\xtilde)^TD_k^{-1}(\zk-P_{X,D_k^{-1}}(\zk))\nonumber\\\nonumber
&&\geq 2\ak (\xk-\xtilde)^T\uk+ 2(\xk-\zk)^TD_k^{-1}(\zk-\xkk) \\\nonumber
&&= 2\ak (\xk-\xtilde)^T\uk + 2(\xk-\zk)^TD_k^{-1}(\zk-\xk)+2 (\xk-\zk)^TD_k^{-1}(\xk-\xkk) \\\nonumber
&&= 2\ak (\xk-\xtilde)^T\uk - 2(\xk-\zk)^TD_k^{-1}(\xk-\zk)+2 \ak(\xk-\xkk)^T\uk \\\nonumber
&&\geq 2\ak (\xk-\xtilde)^T\uk -2 \ak^2\|\uk\|^2\cdot \valeriacorr{\|D_k\|} - 2\ak\|\xk-\xkk\|\cdot\|\uk\|\nonumber\\
&&\geq 2\ak (\xk-\xtilde)^T\uk -2 \valeriacorr{L_k}\ak^2\|\uk\|^2 - 2\ak^2 \valeriacorr{L_k}\|\uk\|^2\nonumber\\
&&\geq 2\ak (\xk-\xtilde)^T\uk -4\ak^2 \valeriacorr{L_k}\rho^2\label{eq4_3}\\
&&\geq 2 \ak(f(\xk)-f(\xtilde)-\ek)-4\ak^2\valeriacorr{L_k}\rho^2\nonumber\\
&&\geq -2 \ak\ek-4\ak^2 \valeriacorr{L_k}\rho^2\nonumber\\
&&\geq -2 L_k\ak\ek-4\ak^2\valeriacorr{L_k}\rho^2\nonumber
\end{eqnarray}
where the first inequality follows from \eqref{ine0}, while we use Lemma \ref{lemma:0} (ii) in the second one, the definition of $\zk$ and the Cauchy--Schwartz inequality in the third one, the assumption \valeriacorr{$\|D_k\|\leq L_k$}
and \eqref{2.2a} in the fourth one, the bound $\|\uk\|\leq \rho$ in the fifth one, the definition of $\epsilon$-subgradient in the sixth one, the fact that $\xtilde\in X^*$ in the seventh one and the inequality $L_k\geq 1$ in the last one.\\
Upon rearranging terms, this yields
\begin{equation}\label{equ1}
\|\xkk-\xtilde\|_{D_k^{-1}}^2\leq \|\xk-\xtilde\|_{D_k^{-1}}^2 + \valeriacorr{5 L_k\ak^2\rho^2} +2 L_k\ak\ek
\end{equation}
Since we have
\begin{equation}\label{new}
\begin{array}{lcl}
\displaystyle\|\xkk-\xtilde\|_{D_k^{-1}}^2&\geq& \displaystyle\lambda_{min}({D_k^{-1}})\|\xkk-\xtilde\|^2\geq \frac 1{L_k} \|\xkk-\xtilde\|^2\\
\displaystyle\|\xk-\xtilde\|_{D_k^{-1}}^2&\leq&\displaystyle\lambda_{max}({D_k^{-1}})\|\xk-\xtilde\|^2\leq L_k \|\xk-\xtilde\|^2
\end{array}
\end{equation}
we can also write
\begin{eqnarray*}
\|\xkk-\xtilde\|^2&\leq& L_k^2\|\xk-\xtilde\|^2 + \valeriacorr{5 L_k^2\ak^2\rho^2} +2 L_k^2\ak\ek\\
&\leq &L_k^2\|\xk-\xtilde\|^2 + \sigma L_k^2\ak^2 +2 L_k^2\ak\ek
\end{eqnarray*}
where \valeriacorr{$\sigma = 5 \rho^2$}. By repeatedly applying the previous inequality we obtain
\begin{eqnarray}
\|\xkk-\xtilde\|^2 &\leq &L_k^2\|\xk-\xtilde\|^2 + \sigma L_k^2\ak^2 +2 L_k^2\ak\ek\nonumber\\
&\leq& L_k^2(L_{k-1}^2\|\x^{(k-1)}-\xtilde\|^2 + \sigma L_{k-1}^2\alpha_{k-1}^2 + 2L_{k-1}^2\alpha_{k-1}\epsilon_{k-1}) + \sigma L_k^2\ak^2 +2 L_k^2\ak\ek\nonumber\\
&\leq& \theta_0^k\|\x^{(0)} -\xtilde\|^2 + \sigma\sum_{j=0}^k\theta_j^k \alpha_j^2 + 2\sum_{j=0}^k\theta_j^k\alpha_j\epsilon_j\label{ine4}
\end{eqnarray}
where $\theta_j^k = \prod_{i=j}^k L_i^2$, $j\leq k$. Since $L_i^2 \geq 1$ we have $1\leq\theta_j^k\leq \theta_{j-1}^k\leq\theta_0^k$, which implies
\begin{eqnarray*}
\|\xkk-\xtilde\|^2
&\leq& \theta_0^k\|\x^{(0)} -\xtilde\|^2 + \sigma\theta_0^k\sum_{j=0}^k \alpha_j^2 + 2\theta_0^k\sum_{j=0}^k\alpha_j\epsilon_j\\
&\leq& M\left(\|\x^{(0)} -\xtilde\|^2 + \sigma\sum_{j=0}^k \alpha_j^2 + 2\sum_{j=0}^k\alpha_j\epsilon_j\right)
\end{eqnarray*}
where the last inequality follows from Lemma \ref{lemma:3}. Thus, by conditions \eqref{hp_sumepsilon}, the sequence $\{\xk\}$ is bounded. In order to show that $\{\xk\}$ converges to a solution of \eqref{minf}, we now consider inequality \eqref{eq4_3}, which, in view of \eqref{new}, results in
\begin{eqnarray*}
\|\xkk-\xtilde\|^2 &\leq &L_k^2\|\xk-\xtilde\|^2 + L_k^8\ak^2\rho^2 +4\ak^2L_k^4\rho^2 -2\ak L_k(\xk-\xtilde)^T\uk\\
&\leq& L_k^2\|\xk-\xtilde\|^2 + \sigma L_k^2 \ak^2 +2\ak L_k(\xtilde-\xk)^T\uk
\end{eqnarray*}
By repeatedly applying the previous inequality we obtain
\begin{eqnarray}
\|\xkk-\xtilde\|^2 &\leq& 
\theta_0^k \|\x^{(0)}-\xtilde\|^2 +\sigma \sum_{j=0}^k\theta_j^k\alpha_k^2 + 2\sum_{j=0}^k\tilde\theta_j^k\alpha_j(\xtilde-x^{(j)})^Tu^{(j)}\label{ine3}
\end{eqnarray}
where $\tilde \theta_j^k = \theta_j^k/L_j$.
Since $\xtilde\in X^*$ and $\u^{(j)}\in\partial_\ek f(\x^{(j)})$, for all $j\geq 0$ we have
\begin{equation}\label{defesubgrad}
f(\x^{(j)})\geq f(\xtilde)\geq f(\x^{(j)})+(\xtilde-\x^{(j)})^T\u^{(j)} -\epsilon_j \end{equation}
Hence, $(\xtilde-\x^{(j)})^T\u^{(j)} \leq \epsilon_j$ for all $j\geq 0$.\\
Now, we show that $(\xtilde-\x^{(j)})^T\u^{(j)}\rightarrow 0$ for $j\rightarrow \infty$.
To this end, assume by contradiction that $(\xtilde-\x^{(j)})^T\u^{(j)} < - \epsilon$ for some $\epsilon > 0$. Then, by inequality \eqref{ine3} we obtain
\begin{eqnarray*}
\|\xkk-\xtilde\|^2
&\leq & M \|\x^{(0)}-\xtilde\|^2 +\sigma M\sum_{j=0}^k\alpha_k^2 - 2\frac {\epsilon} L\sum_{j=0}^k\alpha_j
\end{eqnarray*}
where we applied the inequalities $1\leq\theta_j^k\leq \theta_{j-1}^k\leq\theta_0^k\leq M$, $1\leq L_k\leq L$. Then, taking limits for $k\rightarrow\infty$, by assumption \eqref{hp_aezero} we have an absurdum.
Thus, there exists a subsequence $\{\x^{(k_i)}\}$ such that $\lim_i (\xtilde-\x^{(k_i)})^T\u^{(k_i)} = 0$. Then, from \eqref{defesubgrad} and the assumption \eqref{hp_sumepsilon} it follows that $\lim_i f(\x^{(k_i)}) = f^*$. Since $\{\xk\}$ is bounded, $\{\x^{(k_i)}\}$ is also bounded and, thus it has an accumulation point $\x^{\infty}$. By the continuity of $f(x)$ and by inequality \eqref{defesubgrad} we can conclude that $\x^{\infty}\in X^*$.

Now we show that the whole sequence $\{\xk\}$ converges to $\x^{\infty}$. Let $\delta > 0$; since $\x^{\infty}$ is an accumulation point of $\{\xk\}$ and from \eqref{hp_sumepsilon}, there exists a positive integer $m_\delta$ such that
$\|\x^{\infty}-\x^{(m_\delta)}\|^2\leq \delta/(3M)$, $\sum_{j=m_\delta}^\infty \ak^2< \delta/(3\sigma M)$ and $\sum_{j=m_\delta}^\infty \ak\ek< \delta/(6 M)$. Then, for any $k>m_\delta$, using the same arguments as in \eqref{ine4}, we obtain
\begin{eqnarray*}
\|\xk-\x^{\infty}\|^2
&\leq& \theta_{m_\delta}^{k-1}\|\x^{(m_\delta)} -\x^{\infty}\|^2 + \sigma\sum_{j=m_\delta}^{k-1}\theta_j^{k-1} \alpha_j^2 + 2\sum_{j=m_\delta}^{k-1}\theta_j^{k-1}\alpha_j\epsilon_j\\
&\leq& M\|\x^{(m_\delta)} -\x^{\infty}\|^2 + \sigma M\sum_{j=m_\delta}^{\infty} \alpha_j^2 + 2M\sum_{j=m_\delta}^{\infty}\alpha_j\epsilon_j\\
&\leq& \delta
\end{eqnarray*}
Since $\delta$ can be chosen arbitrarily small, then $\{\xk\}$ converges to $\x^{\infty}$. \endproof
We now discuss some further issues about Theorem \ref{teo:1}, relating our results with the recent literature, in particular with the papers \cite{Combettes-Vu-2013,Combettes-Vu-2014} \silviacorr{and \cite{Nedic-2002}}.\vspace{0.3cm}\\
{\bf Remark 1.} From equation \eqref{equ1}, observing that
\begin{eqnarray*}
\|\xkk-\xtilde\|^2_{D_k^{-1}}&\geq& \lmin(D_k^{-1})\|\xkk-\xtilde\|^2\\
&=& \frac{\lmin(D_k^{-1})}{\lmax(D_{k+1}^{-1})}\lmax(D_{k+1}^{-1})\|\xkk-\xtilde\|^2 \\
&\geq&\lmin(D_k^{-1})\lmin(D_{k+1})\|\xkk-\xtilde\|_{D_{k+1}^{-1}}^2\\
&\geq&\frac{1}{L_kL_{k+1}}\|\xkk-\xtilde\|_{D_{k+1}^{-1}}^2
\end{eqnarray*}
we obtain
$$ \|\xkk-\xtilde\|_{D_{k+1}^{-1}}^2\leq \zeta_k\|\xk-\xtilde\|_{D_k^{-1}}^2 + \xi \zeta_k\ak^2 +2L \zeta_k\ak {\ek}
$$
where $\zeta_k = \sqrt{(1+\gamma_k)(1+\gamma_{k+1})}$ and $\xi=5L\rho^2$. By the assumptions made on $\{\gamma_k\}$, the sequence $\{\zeta_k\}_\kinN$ is bounded. We can also set $\zeta_k= 1+\eta_k$, with $\eta_k= \sqrt{(1+\gamma_k)(1+\gamma_{k+1})}-1$, and observe that the series $\sum \eta_k$ and $\sum \gamma_k$ have the same behaviour, thanks to the limit $\lim_{z\rightarrow 0} (\sqrt{1+z}-1)/z=1/2$. Then, since the assumption \eqref{hp_gammak}, we can conclude that $\sum \eta_k$ is a convergent series.

Thus, the sequence $\{\xk\}_{k\in\mathbb N}$ is {\em quasi-Fej\'{e}r monotone} with respect to $X^*$ relative to $\{D_k^{-1}\}_\kinN$, in the sense of \cite[Definition 3.1]{Combettes-Vu-2013} and we could apply Proposition 3.2 in \cite{Combettes-Vu-2013} (see also \cite{Combettes-Vu-2014}) to obtain that  $\{\|\xk-\xtilde\|_{D_k^{-1}}\}_{k\in \mathbb N}$ converges and, thus, $\{\xk\}_{k\in \mathbb N}$ is bounded.

However, the techniques we used in the proof of Theorem \ref{teo:1} to show the boundedness of $\{\xk\}$ allow to prove the convergence of the whole sequence $\{\xk\}$ to a solution of \eqref{minf} and are employed also in the following section.

\silviacorr{Variable metric was introduced also in \cite[Chapter 5]{Nedic-2002} in the context of subgradient methods for unconstrained problems (i.e. $X=\R^n$). In this case, setting $D_k = B_kB_k^T$, the scaling matrices are assumed to satisfy $\|B_{k+1}^{-1}B_k\|\geq 1$ and $\prod_{k=0}^\infty\|B_{k+1}^{-1}B_k\|^2<\infty$. 
Even if the second condition is verified under the assumptions of Theorem \ref{teo:1}, we observe that the requirement $\|B_{k+1}^{-1}B_k\|\geq 1$ restricts the choice of the scaling matrix, strictly connecting the metrics adopted in two successive iterates.
} \vspace{0.3cm}\\  
\silviacorr{{\bf Remark 2.} One of the crucial assumptions of Theorem \ref{teo:1} concerns the boundedness of the sequence $\{\uk\}$, $\uk\in \partial_{\ek}f(\xk)$. This assumption is satisfied for example when $\diam(\dom(f^*))$ is finite. Indeed, it holds $\dom(f^*) = \bigcup_{\x\in\R^n} \partial_\epsilon f(\x)$, for every $\epsilon>0$ (see Remark 2 in \cite{Villa-etal-2013}). Moreover, if $\diam(\dom(f^*))=M$ for some $M>0$, by definition of $\epsilon$-subdifferential we obtain $|f(\x)-f(\z)|\leq M\|\x-\z\| + \epsilon $. Since $\epsilon$ is arbitrary, it follows that $f$ is Lipschitz continuous with constant $M$.}\vspace{0.3cm}\\
{\bf Remark 3.} The convergence analysis of Theorem \ref{teo:1} applies also to further variants of the iteration \eqref{method}, as for example the following one:
\begin{equation}
\xkk = P_{X,D_k^{-1}}\left(\xk-\frac{\ak}{\max(1, \|\uk\|_{D_k})} D_k {\uk}\right)
\label{method1}
\end{equation}
The convergence of the method \eqref{method1} can be analyzed as follows.
\begin{Cor}\label{cor:0}
Let $\{\xk\}_{k\in \mathbb N}$ a sequence satisfying \eqref{method1} such that the assumptions of Theorem \ref{teo:1} are satisfied. Then, $\{\xk\}_{k\in \mathbb N}$ converges to a solution of problem \eqref{minf}.
\end{Cor}
{\it Proof.} Let us define $\bar{\alpha}_k = \frac{\ak}{\max(1, \|\uk\|_{D_k})}$. From the inequalities
\begin{equation}\nonumber
1\leq \max(1, \|\uk\|_{D_k})\leq \max(L^{\frac{1}{2}}\rho,1)
\end{equation}
it follows that
\begin{eqnarray*}
&&\sum \bar{\alpha}_k= \sum \frac{\ak}{\max(1, \|\uk\|_{D_k})} \geq  \frac{1}{\max(L^{\frac{1}{2}}\rho,1)} \sum \ak =\infty\nonumber \\
&&\sum \bar{\alpha}_k^2=  \sum \frac{\ak^2}{\max(1, \|\uk\|_{D_k})^2} \leq \sum \ak^2 <\infty \nonumber \\
&&\sum \bar{\alpha}_k \ek =  \sum \frac{\ak}{\max(1, \|\uk\|_{D_k})} \ek \leq  \sum \ak \ek <\infty \nonumber \\
\end{eqnarray*}
Then, we can invoke Theorem \ref{teo:1} to obtain the result. \endproof
The results in this Section can be exploited for the practical implementation of the methods \eqref{method} and \eqref{method1}, since they indicates how to choose the sequence $\{\ak\}_\kinN$, and they are employed also in the convergence analysis of \eqref{method1} equipped with an adaptive stepsize rule, as the one described in the following.

\section{Convergence analysis with dynamic stepsize rule}\label{sec3}

A critical point for the implementation of the methods \eqref{method} and \eqref{method1} is how to select the sequence $\{\ak\}$; a practical strategy to obtain good performances is still an open problem, since they are in general quite sensitive to this choice \cite{BonettiniRuggiero12}. \valeriacorr{Borrowing the ideas of \cite{Brannlund95} and \cite{Goffin99}, in this section we describe 
a \emph{level algorithm} that allows to adaptively compute a stepsize $\ak$ of the form \eqref{dynamic} in the iteration \eqref{method1}. In our scheme we introduce the use of the $\epsilon$-subgradient of $f$  at the current iterate (instead of the subgradient) and a variable metric.}
The resulting algorithm is detailed in Algorithm \ref{levelalgorithm}, whose underlying assumption is that, for any given $\epsilon_k\geq0$, we are able to provide an element $\uk$ of the set $\partial_{\epsilon_k}f(\xk)$.
\begin{algorithm}[h]
\caption{Scaled $\epsilon$-Subgradient Level Algorithm (SSL)}\label{levelalgorithm}
\valeriacorr{Choose $B>0$, $\nu_1,\nu_2\in (0,1)$, $\frec_{-1}=\infty$; $k=0$, $l=0$, $k(l)=0$, $\delta_0>0$; choose $x^{(0)}\in X$.
\begin{itemize}
\item[]\textsc{For} $k=0,1,2,...$
\begin{AlgorithmSteps}[4]
\item[1] Computation of $f(\xk)$
\item[2] If $f(\xk)<\frec_{k-1}$, then $\frec_{k}=f(\xk)$ 
else  $\frec_{k}=\frec_{k-1}$ 
\item[3] If  $f(\xk)<\frec_{k(l)}- \nu_1 \delta_l$, then $k(l+1)=k$, $\sigma_k=0$, $\delta_{l+1}=\delta_l$, $l=l+1$ and go to Step 5.
\item[4] If  $\sigma_k>B$, then $k(l+1)=k$, $\sigma_k=0$, $\delta_{l+1}=\nu_2 \delta_l$, $l=l+1$.
\item[5] Set $\flev_k=\frec_{k(l)}-\delta_l$
\item[6] Update the stepsize and compute the new iterate
$$\alpha_k=\frac{f(\xk) -\flev_k}{\max(1,\|\uk\|_{D_k})}$$
    \begin{equation}
    \xkk= P_{X,D_k^{-1}}\left(\xk -\alpha_k D_k\frac{\uk}{\max(1,\|\uk\|_{D_k})}\right) \label{updatexx}
    \end{equation}
\item[7] $\sigma_{k+1}=\sigma_k+\alpha_k $ and go to Step 1.
\end{AlgorithmSteps}
\end{itemize}}
\textsc{End}
\end{algorithm}
\valeriacorr{In Algorithm \ref{levelalgorithm}, we have $\frec_{k}= \min_{i=0,...,k}f(\x^{(i)})$, 
 while $l$ is the number of times that the value $\flev$ has been updated and $k(l)$ is the iteration where the $l$-th updating occurred. Finally, $\sigma_k$ is the cumulative path length between two successive updates of $\flev$.}\\
Steps 2-5 aim to provide in $\flev_k$ an estimate of the optimal function value at the iterate $k$, which is used as target level for the successive iterates until the objective function value is sufficiently close to it or the iterates move through a long path without approaching it. In the first case, i.e. when the inequality at Step 3 is satisfied, $\flev_k$ is reduced at Step 5 by subtracting the positive quantity $\delta_l$ to the best value obtained so far, $\frec$. In the other case, when the inequality at Step 4 is satisfied, the estimated difference from the optimal value $\delta_l$ is reduced and, as a consequence of Step 5, the target level $\flev_k$ is increased.

One of the main differences between the stepsize computed by Algorithm \ref{levelalgorithm} and the square summable sequence considered in the previous section is that the former one does not necessarily converge to zero.

\valeriacorr{In the rest of this section we prove that the sequence $\{f(\xk)\}$, where $\xk$ is computed by Algorithm \ref{levelalgorithm}, converges to the minimum of $f$, using similar techniques as in \cite{Nedic01}. Before giving the main result, whose proof also exploits the results in Section \ref{sec2}, we recall the following technical lemma. We omit the proof, since it runs as that of Lemma 2.2 in \cite{Nedic01}.}
\valeriacorr{\begin{Lemma}\label{lemma3}
Assume that the set $X^*$ of the solutions of (\ref{minf}) is nonempty. Assume that there exists a positive constant $\rho$ such that $\|\uk\|\leq \rho$, $\uk\in \partial_\ek f(\xk)$,  and a sequence of positive numbers $\{L_k\}$ such that $\|D_k\|\leq L_k$, $\|D_k^{-1}\|\leq L_k$, with $1\leq L_k\leq L$ for some positive constant $L$, for all $k\geq 0$.
Given $B>0$ and $\{\epsilon_k\}$ such that 
$\ek\rightarrow 0$ as $k\rightarrow \infty$ in Algorithm \ref{levelalgorithm}, we have ${l}\rightarrow \infty$ and $\delta_l\rightarrow 0$ as $l\rightarrow \infty$.
\end{Lemma}}
%
The following Theorem can be considered as a generalization of Proposition 2.7 in \cite{Nedic01}, which only deals with the case $D_k=I$, $\epsilon_k=0$ for all $k$.

\begin{Thm}\label{theorem3}
Under the same assumptions of the previous Lemma, if \eqref{hp_gammak} holds, for SSL we have $\tilde{f}=\inf_{k\geq 0} f(\xk)= f(x^*)$, with $x^*\in X^*$.
\end{Thm}
\noindent \valeriacorr{{\it Proof.} The first part of the proof aims to show that $\sum_j \alpha_j=\infty$ and runs as Proposition 2.7 in \cite{Nedic01}. For sake of completeness, we report below the detailed derivation of the result.}\\
From Lemma \ref{lemma3}, $\delta_l\rightarrow 0$ as $l\rightarrow \infty$. Let $S$ be given by $S=\{l\in \{1,2,...\}, \delta_l=\nu_2 \delta_{l-1} \}$. Then, from Step 4 and 6 of Algorithm \ref{levelalgorithm}, we obtain
$$\sigma_k=\sigma_{k-1}+\alpha_{k-1} =\sum_{j=k(l)}^{k-1}  \alpha_j $$
so that $k(l+1)=k$ and $l+1\in S$ whenever $\sum_{j=k(l)}^{k-1}  \alpha_j>B $ at Step 4. Hence
$$ \sum_{j=k(l-1)}^{k(l)-1}  \alpha_j >{B}\ \ \ \ \ \forall l\in S $$
and since the cardinality of $S$ is infinite, we have
\begin{equation}\label{new2}
\sum_{k=k(l)}^{\infty}  \alpha_k\geq  \sum_{l\geq \overline{l}, l\in L} \sum_{j=k(l-1)}^{k(l)-1} \alpha_j > \sum_{l\geq \overline{l}, l\in S} {B}=\infty
\end{equation}
Now in order to obtain a contradiction, assume that $\tilde{f}>f(x^*)$, so that for some $\tilde{y}\in X$ and some $\eta>0$
\begin{equation}
\tilde{f}-\eta \geq f(\tilde{y})
\label{absurde}\end{equation}
\valeriacorr{Since $\delta_l\rightarrow 0$ and $\ek\rightarrow 0$, there are large enough $\overline{l}$ and $\overline{k}$ such that, for all $l\geq \overline{l}$ and $k\geq \overline{k}$, we have $\delta_l<\eta/2$ and $\ek<\eta/2$ ; then for all $k\geq \tilde{k}=\max(k(\overline{l}),\overline{k})$
$$\flev_k-\ek=\frec_{k(l)}-\ek-\delta_l> \tilde{f} -\eta \geq f(\tilde{y})$$}
From this inequality, by Lemma \ref{lemma:0} (iii), the definition of $\epsilon$-subgradient, the definition of $\alpha_k$ we obtain
\begin{eqnarray*}
\|x^{(k+1)}-\tilde{y}\|_{D_k^{-1}}^2 &\leq & \|x^{(k)}-\tilde{y}\|_{D_k^{-1}}^2 -2\frac{\alpha_k}{\max(1,\|\uk\|_{D_k})} (f(\xk)-f(\tilde{y}) -\ek)+ \alpha_k^2  \\
&\leq & \|x^{(k)}-\tilde{y}\|_{D_k^{-1}}^2 -2\frac{\alpha_k}{\max(1,\|\uk\|_{D_k})} (f(\xk)-\flev_k )+ \alpha_k^2  \\
&\leq & \|x^{(k)}-\tilde{y}\|_{D_k^{-1}}^2 -  \alpha_k^2
\end{eqnarray*}
In view of \eqref{new} with $\tilde{\x}=\tilde{\y}$ and $L_k\geq 1$, we can write
\begin{equation}
\|x^{(k+1)}-\tilde{y}\|^2 \leq  L_k^2 \|x^{(k)}-\tilde{y}\|^2 -  \alpha_k^2
\label{100}\end{equation}
By repeatedly applying the previous inequality
we obtain
\valeriacorr{$$ \|x^{(k+1)}-\tilde{y}\|^2 \leq \theta_{\tilde{k}}^k \|x^{(\tilde{k})}-\tilde{y}\|^2 -  \sum_{j=\tilde{k}}^k \theta_{j+1}^k \alpha_j^2
$$}
where $\theta_{j}^k=L_j^2\cdot...\cdot L_k^2$; since $\theta_{j}^k\leq \theta_{0}^k\leq  M$, where $M$ is a positive constant (see Lemma \ref{lemma:3}) and $\theta_{j}^k\geq 1$
we have
\valeriacorr{$$ \sum_{\tilde{k}}^\infty \alpha_j^2 \leq M \|x^{(\tilde{k})}-\tilde{y}\|^2$$}
and consequently \valeriacorr{$\sum_{k=\tilde{k}}^\infty \alpha_k^2 <\infty$. Then $\alpha_k\rightarrow 0$ as $k\rightarrow \infty$ and, from \eqref{new2}, $ \sum_{k=\tilde{k}}^{\infty}  \alpha_k =\infty$.} \\
Now we show that $\sum \ak\ek < \infty$.
Indeed, since $\ek\rightarrow 0$ as $k\rightarrow \infty$, there exists $\bar{k}$ such that $2\ek<\eta$ for $k\geq \bar{k}$, where $\eta$ is such that (\ref{absurde}) holds.
We consider the inequality
\begin{equation}
\|x^{(k+1)}-\tilde{y}\|_{D_k^{-1}}^2 \leq  \|x^{(k)}-\tilde{y}\|_{D_k^{-1}}^2 +\alpha_k^2 -2\frac{\alpha_k}{\max(1,\|\uk\|_{D_k})}{\uk}^T(x^{(k)}-\tilde{y})
\label{101}
\end{equation}
For the convexity of $f$, the inequality (\ref{absurde}) and inequality $2\ek<\eta$, we have
$$ f(x^{(k)})+ {\uk}^T (\tilde{y}-x^{(k)})-\ek\leq f(\tilde{y})\leq \tilde{f}-\eta \leq f(x^{(k)})-2\ek$$
Then we have
$$  {\uk}^T (\tilde{y}-x^{(k)})\leq -\ek$$
Using this inequality in  (\ref{101}), we obtain
$$\|x^{(k+1)}-\tilde{y}\|_{D_k^{-1}}^2 \leq  \|x^{(k)}-\tilde{y}\|_{D_k^{-1}}^2 +\alpha_k^2-2\frac{\alpha_k \ek}{\max(1,\|\uk\|_{D_k})} $$
Using the same arguments as above, we obtain
$$\|x^{(k+1)}-\tilde{y}\|^2 \leq L_k^2 \|x^{(k)}-\tilde{y}\|^2 + L_k^2 \alpha_k^2 - 2 \frac{\alpha_k \ek}{\max(1,L^{\frac{1}{2}}\rho)} $$
By repeatedly applying the previous inequality we have
\begin{eqnarray*}
&& \|x^{(k+1)}-\tilde{y}\|^2 \leq \theta_{\bar{k}}^k \|x^{(\bar{k})}-\tilde{y}\|^2 + \theta_{\bar{k}}^k \sum_{j=\bar{k}}^k \alpha_j^2 - \frac{2}{\max(1,L^{\frac{1}{2}}\rho)} \sum_{j=\bar{k}}^k {\alpha_j \epsilon_j}\\
&& \leq M \left(\|x^{(\bar{k})}-\tilde{y}\|^2+ \sum_{j=\bar{k}}^k \alpha_j^2\right) - \frac{2}{\max(1,L^{\frac{1}{2}}\rho)} \sum_{j=\bar{k}}^k {\alpha_j \epsilon_j}
\end{eqnarray*}
Then we have
$$  \sum_{j=\bar{k}}^\infty {\alpha_j \epsilon_j} \leq \frac{M}{2} \max(1,L^{\frac{1}{2}}C) \left(\|x^{(\bar{k})}-\tilde{y}\|^2+ \sum_{j=\bar{k}}^\infty \alpha_j^2\right)<\infty$$
According to Corollary \ref{cor:0} we have $\tilde{f}=f^*$ that contradicts (\ref{absurde}).
\endproof

Further generalizations of Algorithm \ref{levelalgorithm} could be included in the analysis of the previous theorem following \cite[p.122]{Nedic01}, where the authors suggest some modifications of Steps 2, 3 and 4 allowing a variable path bound $B$ and different strategies to update the parameter $\delta_l$. For sake of simplicity we omit here these details.

\section{A Scaled Primal--Dual Hybrid Gradient Method}\label{sec4}

The aim of this Section is to present a concrete example of the method \eqref{method} for the problem \eqref{minf0f1}, which for convenience is reproduced below,
\begin{equation}\label{minsum}
\min_{\x\in X} f(\x) \equiv f_0(\x)+ f_1(A\x)
\end{equation}
where $A\in \mathbb \R^{m\times n}$, $f_0(\x), f_1(\x)$ are convex, proper, lower semicontinuous functions such that $\diam(\dom (f_1^*))$ is finite and $f_1^*(\y)$ is the Fenchel dual of $f_1$. 

We propose the following Scaled Primal--Dual Hybrid  Gradient (SPDHG) method for the solution of \eqref{minsum}
\begin{eqnarray}
\ykk&=& (I+\tau_k \partial f_1^*)^{-1}(\yk + \tau_k A\xk) \label{updatey}\\
\xkk&=& P_{X,D_k^{-1}}(\xk -\ak D_k(\gk + A^T\ykk)) \label{updatex}
\end{eqnarray}
where $\gk\in \partial_{\delta_k}f_0(\xk)$, for some $\delta_k\geq 0$, $\{\tau_k\}$, $\{\ak\}$ are the dual and primal steplength sequences respectively and $(I+\tau_k \partial f_1^*)^{-1}$ is the resolvent operator \cite{Rokafellar70} of $f_1^*$ defined as
$$(I+\tau_k\partial f_1^*)^{-1}(\x) = \mbox{arg}\min_\z f_1^*(\z) +\frac 1{2\tau_k}\|\z-\x\|^2. $$
The basic result allowing to consider the method \eqref{updatey}--\eqref{updatex} as a special case of a scaled $\epsilon$- subgradient method \eqref{method} is the following Lemma (see \cite[Lemma 1]{BonettiniRuggiero12}).
\begin{Lemma}\label{lemma:1}
Let $\ykk$ defined as in \eqref{updatey}. Then, \silviacorr{$\ykk\in \dom(f_1^*)$ and, thus,} $A^T\ykk\in\partial_{\sigma_k} (f_1\circ A)(\xk)$, where $\sigma_k = f_1(A\xk) + f_1^*(\ykk)-{\ykk}^T A\xk$. Moreover, if there exists a positive number $D$ such that $\diam(\dom (f_1^*))\leq D$, then $\sigma_k\leq (2\tau_k)^{-1}D^2$.
\end{Lemma}
Thus, recalling the additivity of the $\epsilon$-subgradient, we can conclude that $\uk = \gk + A^T\ykk\in \partial_{\ek}f(\xk)$, where $\ek = \delta_k +\sigma_k$.

Motivated by the previous observation, building on the material developed in Sections \ref{sec2} and \ref{sec3}, we discuss two stepsize selection strategies for the method \eqref{updatey}--\eqref{updatex}, providing two different SPDHG implementations.
In the first case, we assume that $\{\tau_k\}_\kinN$, $\{\ak\}_\kinN$, $\{L_k\}_\kinN$ are prefixed sequences.
The following Corollary, which is a consequence of Lemma \ref{lemma:1} and of Theorem \ref{teo:1}, gives conditions on these sequences ensuring the convergence of the sequence $\{\xk\}$ generated by SPDHG to a point $\x^*\in X^*$.
\begin{Cor}\label{cor:1}
Let $\{\xk\}$ be the sequence generated by iteration \eqref{updatey}-\eqref{updatex}. Assume that $\gk\in \partial_{\delta_k} f_0(\xk)$ and that there exists $\rho>0$ such that $\|\gk\|\leq\rho$ for all $k$. Let the steplength sequences $\{\tau_k\}_\kinN$, $\{\ak\}_\kinN$ and the scaling matrix bounds $\{L_k\}_\kinN$ satisfy
\begin{equation}\label{rev1}
\ak = {\mathcal O}\left(\frac 1 {k^p}\right),\ \ \ \tau_k = {\mathcal O} (k^p), \ \ \ L_k = \sqrt{1 + {\mathcal O}\left(\frac 1 {k^q}\right)} \ \ \ \frac{1}{2}<p\leq 1,\ \ q > 1.
\end{equation}
Moreover, assume that $\delta_k$ converges to zero at least as $\frac{1}{\tau_k}$. If the set of the solutions of \eqref{minf} is nonempty and $\diam(\dom(f_1^*))$ is finite, then, the sequence $\{\xk\}$ converges to a solution of \eqref{minf}.
\end{Cor}
\silviacorr{\proof As observed above, we have $\uk = \gk + A^T\ykk\in \partial_{\ek}f(\xk)$, where $\ek = \delta_k +\sigma_k$. Since $\diam(\dom(f_1^*))$ is finite, we can apply Lemma \ref{lemma:1} obtaining $\sigma_k\leq (2\tau_k)^{-1}D^2$. By the assumption \eqref{rev1} on $\tau_k$ and on $\delta_k$ we obtain that  $\epsilon_k= {\mathcal O}\left(\frac 1 {k^p}\right)$ and, as a consequence, $\ak\epsilon_k = {\mathcal O}\left(\frac 1 {k^{2p}}\right)$. Since $\frac{1}{2}<p\leq 1$ and $q > 1$, all assumptions \eqref{hp_ezerovera}--\eqref{hp_gammak} of Theorem \ref{teo:1} are satisfied and we obtain the result.\endproof}
On the other side, the SSL procedure for dynamically computing the primal stepsize $\ak$ can be also implemented.
\valeriacorr{In this case, we have to provide a sequence $\{L_k= \sqrt{1+\gamma_k}\}_\kinN$ such that $\sum \gamma_k< \infty$.
} For sake of simplicity we assume $\delta_k = 0$.

\valeriacorr{In this case the value $\ek= f_1(A\xk) +f_1^*(y^{(k+1)}) - {y^{(k+1)}}^T A \xk$ is controlled by the dual stepsize $\tau_k$ (see Lemma \ref{lemma:1}).}
Then
it is possible to compute $\alpha_k$ by Steps 2-5 of Algorithm \ref{levelalgorithm} and the next iterate $\xkk$.\\
Finally, Theorem \ref{theorem3} allows to conclude that the sequence $\{f(\xk)\}_\kinN$ generated by SPDHG method combined with SSL algorithm converges to the minimum of $f(\x)$ in \eqref{minsum}, as stated in the following Corollary under the same boundedness assumptions as in Corollary \ref{cor:1}.
\valeriacorr{\begin{Cor}\label{cor:2}
Let $\{\xk\}_\kinN$ be the sequence generated by Algorithm \ref{levelalgorithm}.
Here $\uk=\gk + A^T\ykk$ in \eqref{updatexx}, $\gk\in \partial f_0(\xk)$ and $\ykk$ is computed as in \eqref{updatey}. Assume that $\lim_{k\rightarrow \infty}\tau_k= \infty$, $L_k = \sqrt{1 + {\mathcal O}\left(\frac 1 {k^q}\right)}$, with $q>1$, and that there exists $\rho>0$ such that $\|\gk\|\leq\rho$. If the set of the solutions of \eqref{minf} is nonempty and $\diam(\dom(f_1^*))$ is finite, then the sequence $\{f(\xk)\}_\kinN$ converges to $f(\x^*)$, with $\x^*\in X^*$.
\end{Cor}}
\silviacorr{\proof Since $\gk\in \partial f_0(\xk)$, by Lemma \ref{lemma:1} we have $\uk\in \partial_{\ek}f(\xk)$, where $\ek= f_1(A\xk) +f_1^*(y^{(k+1)}) - {y^{(k+1)}}^T A \xk$. Since $\diam(\dom(f_1^*))$ is finite, we can apply the second part of Lemma \ref{lemma:1} obtaining $\ek\leq (2\tau_k)^{-1}D^2$, for a positive constant $D$ such that $\diam(\dom(f_1^*))\leq D$.
Since
$\lim_{k\rightarrow \infty}\tau_k= \infty$, we have $\lim_{k\rightarrow \infty}\epsilon_k= 0$ and by Theorem \ref{theorem3} we obtain the result.\endproof}

\section{Application: edge preserving deblurring of Poisson images}\label{sec5}

In this section we further specialize the SPDHG method, by focusing on a specific application in the image restoration context. Our aim is to suggest a strategy to compute a suitable scaling matrix $D_k$, fully defining the algorithm; as observed by several authors, this choice should be driven according to the specific problem features, such as the structure of the constraints and objective function \cite{Bonettini2009,Lanteri2002,Zanella2009}.

For these reasons, we describe first some details of the image reconstruction problems which, in the Bayesian framework, can be formulated as a constrained convex minimization problems of the form \eqref{minsum}.
For these problems, the function $f_0(\x)$ measures the data discrepancy and should be chosen according to the noise statistics: in particular, when the data suffer from Poisson noise, the Maximum Likelihood principle leads to the generalized Kullback--Leibler divergence
\begin{equation}\label{KL}
f_0(\x) = \sum_{i=1}^n g_i\log\frac{g_i}{(H\x)_i+b}+ (H\x)_i +b - g_i \end{equation}
where $ g\in \mathbb R^n$ is the observed image, $H\in \mathbb R^{n\times n}$ represents the blurring operator while $b\in\R$ is a nonnegative background term. Standard assumptions on $H$ are that it has nonnegative entries and $H^Te> 0$, where $e\in\R^n$ is the vector of all ones.

Since the entries of the unknown vector $\x$ represents the image pixels, a meaningful solution is obtained by defining the constraint set as the nonnegative orthant, i.e. $X=\{\x\in\R^n:x_i\geq0\} $.\\
On the other side, $f_1(A\x)$ plays the role of a regularization term enforcing suitable properties on the solution of \eqref{minsum}. Typically, to preserve the edges in the solutions of \eqref{minsum}, $f_1(A\x) $ can be chosen as
\begin{equation}\label{TV}
f_1(A\x) = \beta TV(x),\ \ \ TV(x)=\sum_{i=1}^n \|A_i\x\|, \ \ \ A_i\in\mathbb R^{2\times n}
\end{equation}
where $TV(x)$ is the discrete, nonsmooth, Total Variation (TV) functional, $\beta$ is a positive regularization parameter and $A_i\in \mathbb \R^{2\times n}$ is defined such that $A_i\x$ represents the discrete gradient of the image $\x$ at the pixel $i$. In these settings, the matrix $A$ is defined by blocks as $A= \begin{pmatrix}A_1^T& A_2^T& \cdots& A_n^T\end{pmatrix}^T\in \mathbb \R^{2n\times n}$.\\
In order to simplify the notation, we assume that $\x\in \R^n$ is a $N\times N$ image, i.e. $n=N^2$ and we will indicate the component $x_\ell$, $\ell = 1,\cdots, n$, also as $x_{i,j}$, $i,j=1,\cdots,N$, with the correspondence $j =\lfloor(\ell-1)/N\rfloor + 1$, $i = \ell- \lfloor(\ell-1)/N\rfloor\cdot N$, where $\lfloor\cdot\rfloor$ denotes the integer quotient. With this notation, the $\ell$--th discrete gradient of the image $\x$ can be written as
$$A_\ell\x = \begin{pmatrix}x_{i+1,j}-x_{i,j}\\ x_{i,j+1}-x_{i,j} \end{pmatrix}
$$
where some boundary conditions are assumed.

In this case, since $f_0$ is differentiable, we define $w^{(k)}=\nabla f_0(\xk)$ in \eqref{updatex}, so that $\delta_k = 0$ for all $k$ in Corollary \ref{cor:1}.
Moreover, the resolvent operator in \eqref{updatey} consists in a simple projection onto the set $ B\times B\times\cdots\times B\subset \R^{2n} $, where $B =\{z\in \R^2: \|z\|\leq 1\} $.

In order to devise a suitable scaling matrix $D_k$ for SPDHG, we adapt to our case the split gradient strategy proposed in \cite{Bertero2008,Lanteri2002} for nonnegatively constrained differentiable problems, which demonstrated to be very effective in several applications \cite{Bertero-etal-2009,Bonettini2010a,Prato2013,Theys-etal-2009,Zanella2009}.

The key point of this approach consists in finding a subgradient decomposition of the form $\uk = V(\xk)-U(\xk)$ with $V(\xk)>0$ and $U(\xk)\geq 0$ for all $k$ and then defining $D_k$ in \eqref{updatex} as a diagonal scaling matrix whose entries are the projection of $\xk_i/V_i(\xk)$ onto the set $[1/L_k,L_k]$.

This strategy has the advantage to agree with the nonnegativity constraints and strongly depends on the form of the subgradient $\uk$.

For a practical implementation of this strategy, we have to find a decomposition of the vector $\uk=\nabla f_0(\xk)+\beta A^T\ykk$ as the difference of two nonnegative terms.

As concerns the first term, the gradient of $f_0$ has the natural decomposition $\nabla f_0(x) = H^Te - H^Tv(x)$, where $v(x)$ denotes the vector with entries $v_i(x)=g_i/(H\x+b)_i$; by the assumptions on $H$, we have $H^Te>0$ and $H^Tv(x)\geq 0$ for all $x\geq 0$.

Thus, it remains to find a decomposition of the vector $A^T\ykk$ in \eqref{updatex}. To this end, we compute the explicit expression of it as a function of $\x^{(j)}$, $j=0,...,k$. We first observe that, if the dual variable is partitioned as
$$\y=\begin{pmatrix} \y_1\\ \y_2\\ \vdots \\ \y_n\end{pmatrix},\ \ \ \y_i\in \R^2, $$
the updating rule \eqref{updatey} can be written as
\begin{eqnarray*}
\tildeyk &=& \yk + \tau_k\beta A\xk\\
\ykk&=& S_k\tildeyk
\end{eqnarray*}
where $S_k$ is a diagonal $2n\times 2n$ matrix with the following diagonal entries
\begin{equation}\label{Sk}
(S_k)_{2i-1,2i-1}=(S_k)_{2i,2i}=\frac{1}{\max\{1,\|\tildeyk_i\|\}}, \ \ \ i=1,...,n
\end{equation}
If the method is initialized with $\y^{(0)} = 0$, the dual variable can be written as
\begin{eqnarray*}
\y^{(0)} &=& 0\\
\y^{(1)} &=& \beta\tau_0S_0A\x^{(0)}\\
\y^{(2)} &=& \beta S_1(\tau_0S_0A\x^{(0)}+\tau_1A\x^{(1)})\\
\y^{(3)} &=& \beta S_2(\tau_0S_1S_0A\x^{(0)}+\tau_1S_1A\x^{(1)}+\tau_2A\x^{(2)})\\
&\vdots&\\
\ykk &=& \beta \sum_{j=0}^k \tau_j\tilde S_j^k A\x^{(j)}
\end{eqnarray*}
where
$$\tilde S_j^k = \prod_{i=j}^k S_i $$
As a consequence, the $\epsilon$-subgradient of $f_1\circ A$ employed in \eqref{updatex} can be expressed as
\begin{equation}\label{rel1}
\beta A^T\ykk = \beta^2\sum_{j=0}^k \tau_j A^T\tilde S^k_jA\x^{(j)}\end{equation}
The following simple lemma, which directly follows from the definition of $A$, indicates a possible decomposition of each term in the summation at the right hand side of \eqref{rel1} as the difference between a positive and a nonnegative term.
\begin{Lemma} \label{lemma:4}
Every matrix--vector product of the form $A^TSA\x$ where $S$ is a $2n\times 2n$ diagonal matrix with positive entries such that $S_{2\ell,2\ell}=S_{2\ell-1,2\ell-1}= s_\ell$, $\ell = 1,\cdots,n$, $\x\geq 0$, can be decomposed as
$$A^TSA\x = V_S\x-U_S\x $$
where
\begin{eqnarray*}
(V_S\x)_{i,j} &=& (2s_{i,j}+s_{i,j-1}+s_{i-1,j})x_{i,j}\geq 0\\
(U_S\x)_{i,j} &=& s_{i,j}(x_{i+1,j}+x_{i,j+1}) +s_{i,j-1}x_{i,j-1}+s_{i-1,j}x_{i-1,j}\geq 0
\end{eqnarray*}
\end{Lemma}
with the correspondence $s_\ell\equiv s_{i,j}$, $j =\lfloor(\ell-1)/N\rfloor + 1$, $i = \ell- \lfloor(\ell-1)/N\rfloor\cdot N$.\\\vspace{0.5cm}\noindent
The $\epsilon$-subgradient of $f$ in \eqref{updatex} can be decomposed as
\begin{equation}\nonumber
\uk = \nabla f_0(\xk) + \beta A^T\ykk= V(\xk)-U(\xk)
\end{equation}
where
\begin{equation}\label{Vxk0}V(\xk) = H^Te +\beta^2 \sum_{j=0}^k\tau_j V_{\tilde S_j^k}\x^{(j)}.\end{equation}
Even if it seems quite complicate, the term
$$V^R(\xk)=\beta^2\sum_{j=0}^k\tau_j V_{\tilde S_j^k}\x^{(j)}$$
can be easily computed in a recursive way, by introducing three auxiliary vectors, as described in Algorithm \ref{GPM}. We also notice that, since the scaling matrix $D_k$ is diagonal and the constraint set $X$ is the nonnegative orthant, the projection $P_{X,D_k^{-1}}(\cdot)$ reduces to the usual Euclidean projection $P_{\geq 0}(\cdot)$ (see Step 7).\\
\begin{algorithm}[t]
\caption{Scaled Primal--Dual Hybrid Gradient (SPDHG)}
\label{GPM}
Choose the starting point $\ve x^{(0)}\in X$ and set $\y^{(0)} = 0$, $\ve p^{(-1)}= \ve q^{(-1)} = \ve r^{(-1)}=0$. Choose the sequences $\{\ak\}_\kinN$, $\{\tau_k\}_\kinN$, $\{\gamma_k\}_\kinN$.
\\[.2cm]
{\textsc{For}} $k=0,1,2,...$ \textsc{do the following steps:}
\begin{itemize}
\item[]
\begin{AlgorithmSteps}[4]
\item[1] Compute $\tildeyk = \yk + \beta \tau_k A\xk$;
\item[2] Compute $s^{(k)}_\ell = \frac{1} {\max\{1,\|\tildeyk_\ell\|\}}$, $\ell=1,...,n$ and define $S_k$ as in \eqref{Sk};
\item[3] Dual update: $\ykk = S_k \tildeyk$;
\item[4] Auxiliary vectors update for the decomposition:
\begin{eqnarray}
p^{(k)}_{i,j} &=& (p^{(k-1)}_{i,j}+\beta^2\tau_k x^{(k)}_{i,j})s_{i,j}^{(k)}\label{updatep}\\
q^{(k)}_{i,j} &=& (q^{(k-1)}_{i,j}+\beta^2\tau_k x^{(k)}_{i,j})s_{i-1,j}^{(k)}\label{updateq}\\
r^{(k)}_{i,j} &=& (r^{(k-1)}_{i,j}+\beta^2\tau_k x^{(k)}_{i,j})s_{i,j-1}^{(k)}\label{updater}\\
& & \mbox{ for } i,j=1,...,N\nonumber
\end{eqnarray}
\item[5] Compute the positive part of the decomposition:
\begin{eqnarray}
V(\xk) &=& H^Te +  (2\ve p^{(k)} + \ve q^{(k)} +\ve r^{(k)})
\label{Vxk}
\end{eqnarray}
\item[6] Compute the scaling matrix:
\begin{eqnarray*}
L_k &=& \sqrt{1+\gamma_k}\\
(D_k)_{\ell,\ell}&=& \min\left\{L_k,\max\left\{L_k^{-1},\frac{x^{(k)}_\ell}{V(\xk)_\ell}\right\}\right\}
\end{eqnarray*}
\item[7] Primal update: $\xkk = P_{\geq 0}(\xk -\ak D_k(\nabla f_0(\xk)+\beta A^T\ykk))$.
\end{AlgorithmSteps}
\end{itemize}
\textsc{End}
\end{algorithm}
By induction it can be shown that the computation of $V(\xk)$ in \eqref{Vxk} actually gives \eqref{Vxk0}. For sake of simplicity, we limit ourselves to show that this is true for $k=0,1$. Indeed, from Lemma \ref{lemma:4} and from \eqref{updatep}--\eqref{updater} we have
\begin{eqnarray*}
V^R(\x^{(0)})_{i,j} &=& \beta^2\tau_0 (V_{S_0}\x^{(0)})_{i,j}\\
&=& \beta^2\tau_0 (2s^{(0)}_{i,j} + s^{(0)}_{i-1,j}+s^{(0)}_{i,j-1}) \x^{(0)}_{i,j}\\
&=& 2p^{(0)}_{i,j} + q_{i,j}^{(0)} + r_{i,j}^{(0)}
\end{eqnarray*}
\begin{eqnarray*}
V^R(\x^{(1)})_{i,j} &=& \beta^2\tau_0 (V_{S_0S_1} \x^{(0)})_{i,j}+ \beta^2\tau_1(V_{S_1} \x^{(1)})_{i,j}\\
&=& \beta^2\tau_0 (2s^{(0)}_{i,j}s^{(1)}_{i,j} + s^{(0)}_{i-1,j}s^{(1)}_{i-1,j}+s^{(0)}_{i,j-1}s^{(1)}_{i,j-1}) \x^{(0)}_{i,j}+\\
& & + \beta^2 \tau_1 (2s^{(1)}_{i,j} + s^{(1)}_{i-1,j}+s^{(1)}_{i,j-1}) \x^{(1)}_{i,j}\\
&=& 2(\beta^2 \tau_0s_{i,j}^{(0)} x^{(0)}_{i,j}+\beta^2 \tau_1x^{(1)}_{i,j})s_{i,j}^{(1)} +\\
& & + (\beta^2 \tau_0s_{i-1,j}^{(0)} x^{(0)}_{i,j}+\beta^2 \tau_1x^{(1)}_{i,j})s_{i-1,j}^{(1)} +\\
& & + (\beta^2 \tau_0s_{i,j-1}^{(0)} x^{(0)}_{i,j}+\beta^2 \tau_1x^{(1)}_{i,j})s_{i,j-1}^{(1)} \\
&=& 2 (p_{i,j}^{(0)} +\beta^2\tau_1x_{i,j}^{(1)})s_{i,j}^{(1)} +(q_{i,j}^{(0)} +\beta^2\tau_1x_{i,j}^{(1)})s_{i-1,j}^{(1)} + (r_{i,j}^{(0)} +\beta^2\tau_1x_{i,j}^{(1)})s_{i,j-1}^{(1)}\\
&=& 2p_{i,j}^{(1)} + q_{i,j}^{(1)} + r_{i,j}^{(1)}
\end{eqnarray*}
Algorithm \ref{GPM} can be adapted for both the stepsize selection strategies described in Sections \ref{sec2} and \ref{sec3}. In the first case, three prefixed sequences $\{\alpha_k\}_\kinN$, $\{\tau_k\}_\kinN$, and $\{\gamma_k\}_\kinN$ satisfying the assumptions of Corollary \ref{cor:1} have to be provided.

In the other case, the SSL procedure for dynamically compute the primal stepsize $\ak$ can be included in Algorithm \ref{GPM}.
Here, only the sequences $\{\tau_k\}_\kinN$, and $\{\gamma_k\}_\kinN$ should be given such that $\lim_{k\rightarrow\infty}\tau_k=\infty$ and $\sum \gamma_k<\infty$.
\section{Numerical experience}\label{sec6}
The aim of our numerical experience is twofold: first, we are interested in evaluating the effect of the scaling on the convergence behaviour of the $\epsilon$-subgradient method. Secondly, we compare the two steplength selection strategies presented in Sections \ref{sec2}-\ref{sec3}.

To this end we consider four different versions of the method \eqref{updatey}--\eqref{updatex}:
\begin{description}
\item[PDHG] corresponds to the choices $D_k=I$ and $\ak$ chosen as an \emph{a priori} diminishing, divergent series, square summable sequence in \eqref{updatex}. It actually consists in the method in \cite{BonettiniRuggiero12};
\item[SPDHG] is Algorithm \ref{GPM} with $\ak$ chosen as an \emph{a priori} diminishing, divergent series, square summable sequence;
\item[SL] is the $\epsilon$-subgradient level method given in Algorithm \ref{levelalgorithm} with $D_k=I$ and  $\uk=\nabla f_0(\xk)+\beta A^T\ykk$, where $\ykk$ is updated as in \eqref{updatey};
\item[SSL] is the same as above but with the scaling matrix $D_k$ defined as at the Step 6 of Algorithm \ref{GPM}.
\end{description}
The numerical experiments described in this section have been performed in
the MATLAB environment (R2012b) on a PC equipped with an Intel Core i7-3517U processor 1.9 GHz,
8 GB RAM.

In our experiments, we consider problem \eqref{minsum} where $f_0,f_1$ are defined in \eqref{KL}-\eqref{TV} and $H$ represents the convolution operator with a given point spread function (psf). Thus, assuming periodic boundary conditions, the matrix--vector products involving $H$ can be computed by the Fast Fourier Transform (FFT) and, by a simple normalization of the psf, we also have $H^Te = He = e$; moreover, in our experiments, $H$ is nonsingular, although very ill conditioned, and $g>0$, so that problem \eqref{minsum} has a unique solution \cite{Figueiredo2010}.

We consider a set of three test problems generated as in \cite{Setzer10} where the data $g$ is obtained with the following procedure: the selected original image $\bar x$ is rescaled so that the maximum pixel intensity is a specified value $I_{max}$. Then, the rescaled image is convolved with the psf and the background $b$ is added. Finally, Poisson noise is introduced by the Matlab \verb"imnoise" function and the simulated data $g$ is obtained after scaling back again by $I_{max}$.

For each test problem, the regularization parameter $\beta$ has been empirically selected by computing the solution of \eqref{minsum} for different values of $\beta$ and choosing that for which we observed the minimum $l_2$ relative distance with respect to  $\bar{x}$.

The features of each test problem are specified below.
\begin{description}
\item{\textbf{\emph{cameraman}}}: the $256\times 256$ original image is the `cameraman' available in the Matlab package, while the psf is a Gaussian function, with standard deviation $1.3$, truncated at the $9\times 9$ central pixels. The other parameters are $I_{max}=1000$, $b=0$, $\beta=0.005$; the $l_2$ relative distance between $\bar x$ and $g$ is $0.1209$, while $g_i\in [4, 250]$.
\item \textbf{\emph{micro}}: the original image is the confocal microscopy phantom of size $128\times 128$ described in \cite{Willett03}, scaled by $10$; the psf is the one in \cite{Willett03} truncated at the $9\times 9$ central pixels. Here we set $I_{max}=1$, $b=0$, {$\beta=0.0477$}; the original image pixels are in the range {$[10, 690]$, the $l_2$ relative distance between $\bar x$ and $g$ is $0.1442$, while $g_i\in[1, 778]$}.
\item \textbf{\emph{phantom}}:  the original image $\bar x$ is the $256\times 256$ Shepp-Logan phantom, generated by the Matlab function \verb"phantom", scaled by a factor 1000, while the psf is a Gaussian function, with standard deviation $3$, truncated at the $9\times 9$ central pixels. In this case we set $I_{max}=1$, {$b=10$  and $\beta= 0.00526$}. The values of the original image are in the
range $[0, 1000]$, the $l_2$ relative distance between $\bar x$ and $g$ is {$0.4643$, while $g_i\in [1, 934]$}.
\end{description}
\noindent

For all test problems we compute the solution $x^*$ of the minimization problem (\ref{minsum}) by running 50000 iterations of  PIDSplit method \cite{Setzer10}. Then, we evaluate the progress toward this solution at each iteration in terms of the $l_2$ relative error from the minimum point and the relative difference from the optimal value
$$e^k=\frac{\|\xk-\ve{x}^*\|}{\|\ve{x}^*\|}\ \ \ f^k=\frac{f(\xk)-f(\x^*)}{f(\x^*)}.$$
Following the assumptions of Corollaries \ref{cor:1} and \ref{cor:2}, we choose the sequences of parameters as follows
$$\tau_k = t_1+t_2k \ \ \ak = \frac{1}{t_3+t_4k}\ \ \gamma_k =\frac{t_5}{k^{1+t_6}} $$
In order to illustrate the effectiveness of the methods, the values $t_i$ have been manually optimized for each test problem to obtain a faster decrease of $e^k$ (see Table \ref{table:1}).\\ Moreover, for the initialization of both SL and SSL, we adopt the rule $$\delta_0 = 0.9\ f(\x^{(0)})$$
while the other parameters are $\nu_1 = \nu_2 = 0.5$, $B=0.9\|\u^{(0)}\|\|D_0\|_{\infty}^{\frac 1 2}$.
\\\noindent
\begin{table}
{\begin{tabular}{ r l l l }
                  & \multicolumn{3}{c}{PDHG}\\ \cline{2-4}
                  & $\tau_k$     &     $\ak$                &   $\gamma_k$\\ \cline{2-4}
\emph{cameraman}  & $0.9+10^{-2}k$       & $(0.04+10^{-5}k)^{-1}$&   - \\
\emph{micro}      & $0.9+10^{-3}k$       & $(0.04+10^{-4}k)^{-1}$ &   -\\
\emph{phantom}    & $0.9+ 10^{-3}k$ & $(0.2+10^{-5}k)^{-1}$ &   -\\\hline
                  & \multicolumn{3}{c}{SPDHG}\\\cline{2-4}
                  & $\tau_k$             &     $\ak$                     &   $\gamma_k$\\ \cline{2-4}
\emph{cameraman}  & $0.5+5\cdot10^{-3}k$ & $(0.5+5\cdot10^{-5}k)^{-1}$   & $10^{13}k^{-2}$\\ 
\emph{micro}      & $0.4+10^{-5}k$       & $(0.4+10^{-5}k)^{-1}$         & $10^{13}k^{-2}$\\ 
\emph{phantom}    & $0.5+10^{-4}k$ & $(0.5+10^{-5}k)^{-1}$         & $10^{13}k^{-2}$\\\hline
\end{tabular}\ \
\begin{tabular}{ll}
\multicolumn{2}{c}{SL}\\ \cline{1-2}
$\tau_k$     & $\gamma_k$\\\hline
$0.5+5\cdot10^{-2}k$   &     -     \\
$0.9+10^{-1}k$  &     -     \\
$0.9+ 10^{-2}k$  &     -     \\ \hline
\multicolumn{2}{c}{SSL}\\ \cline{1-2}
$\tau_k$     & $\gamma_k$\\\hline
$0.7+5\cdot10^{-2}k$   &$10^{13}k^{-2}$\\
$0.9+10^{-2}k$ &$10^{13}k^{-2}$\\
$0.9+10^{-2}k$  &$10^{13}k^{-2}$\\\hline
\end{tabular}}
\caption{Parameter settings.}\label{table:1}
\end{table}
\newcommand*{\factor}{0.29}
\begin{figure}
\begin{center}
\begin{tabular}{ccc}
\includegraphics[scale=\factor]{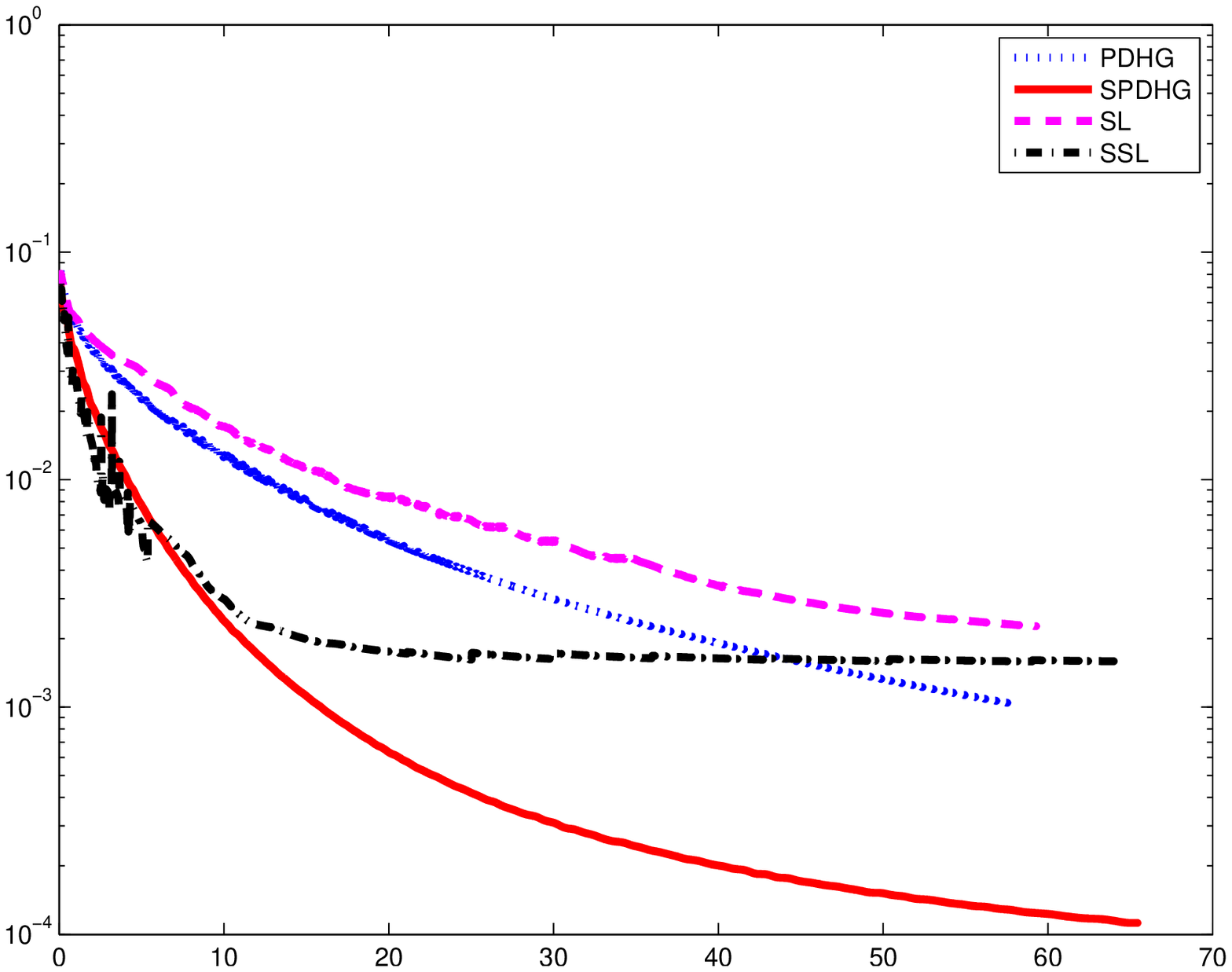}& \includegraphics[scale=\factor]{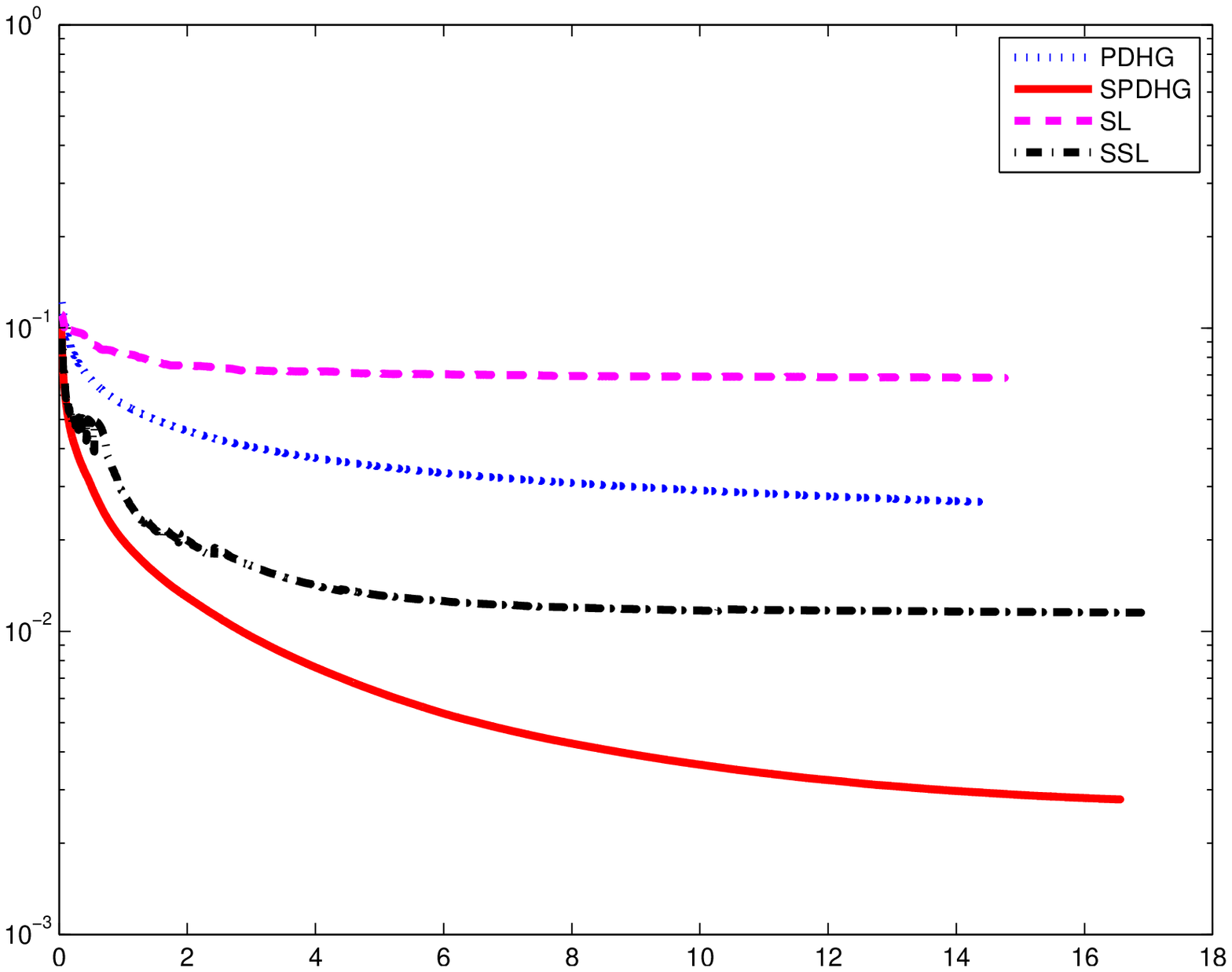}&\includegraphics[scale=\factor]{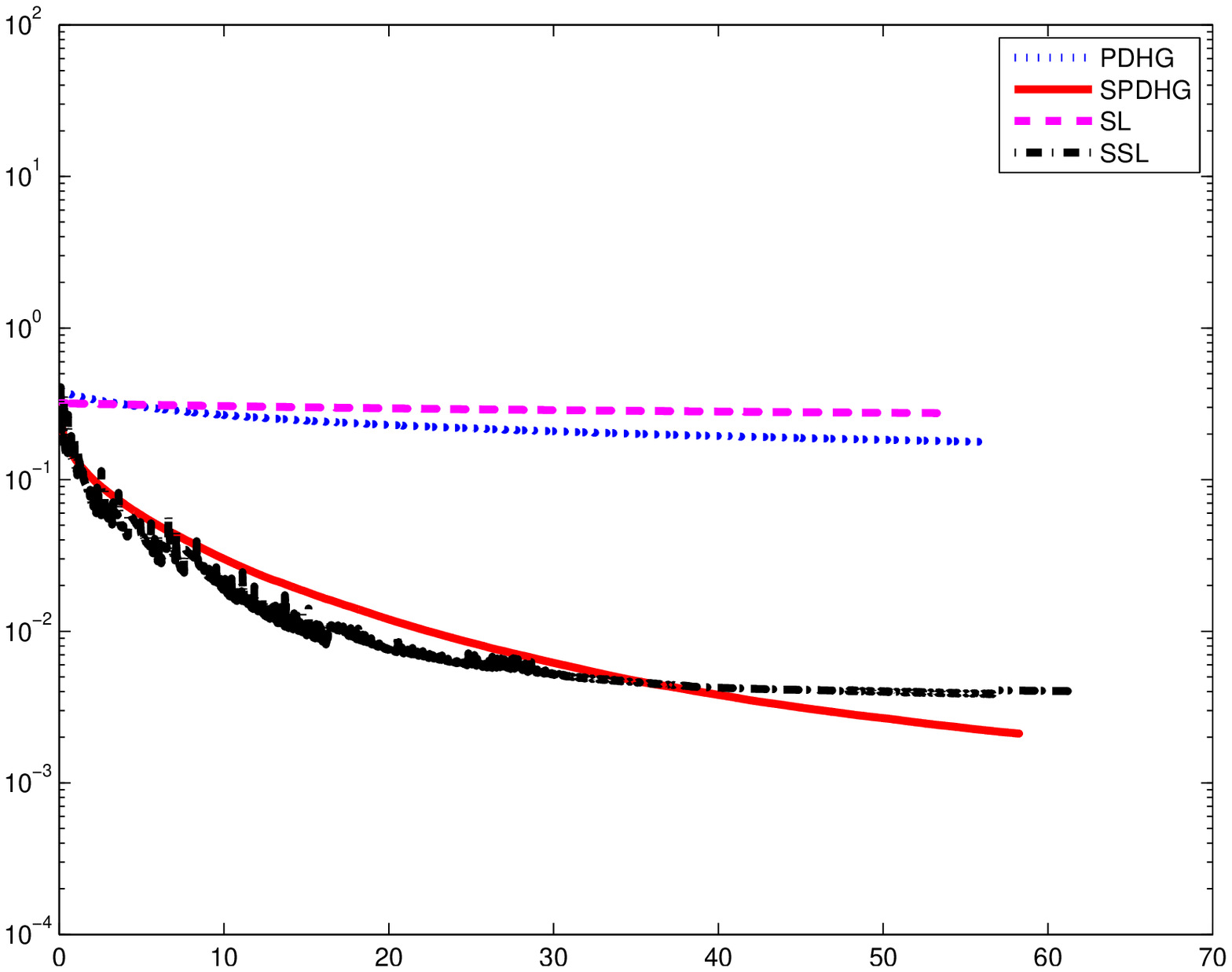}\\
\includegraphics[scale=\factor]{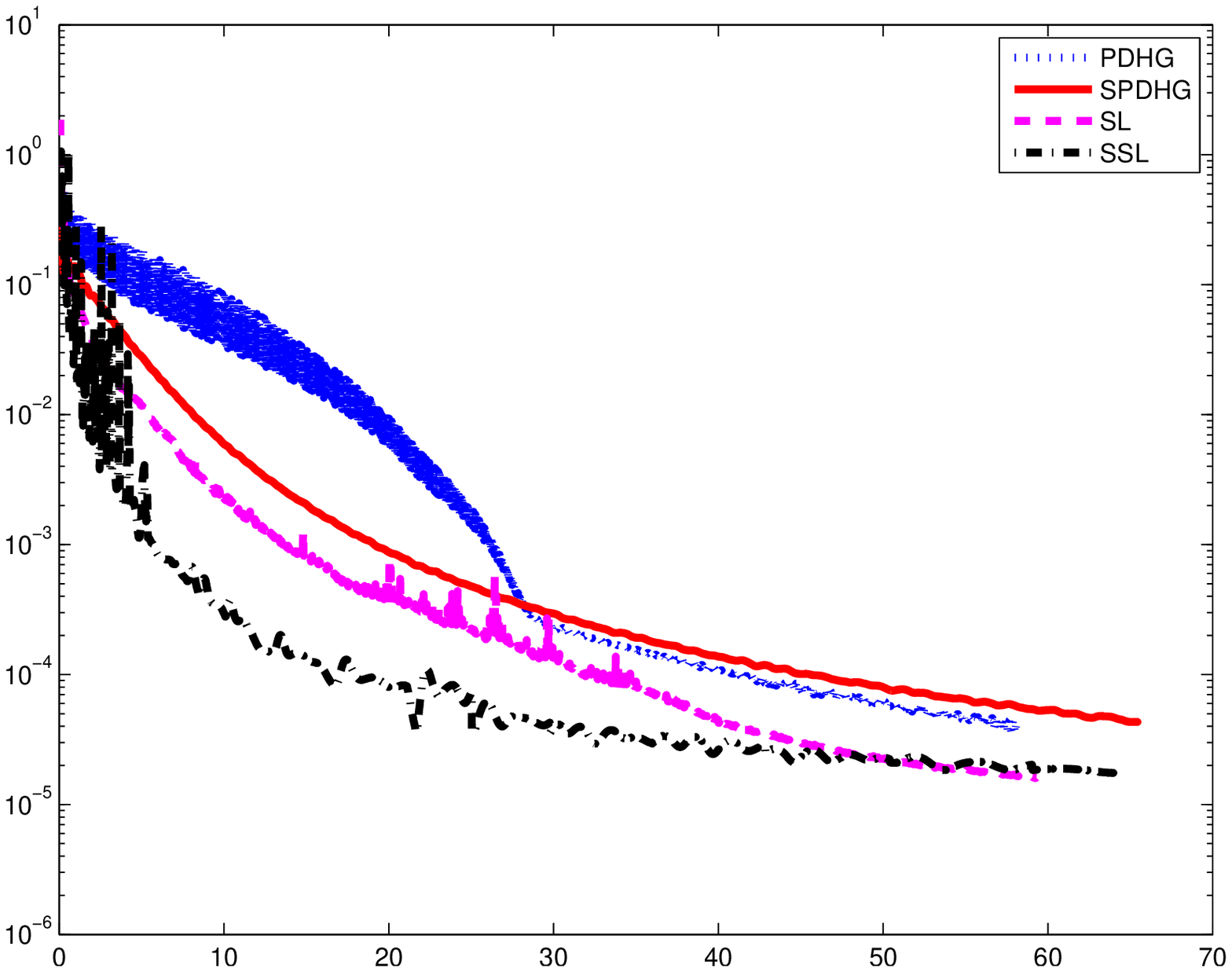}&
\includegraphics[scale=\factor]{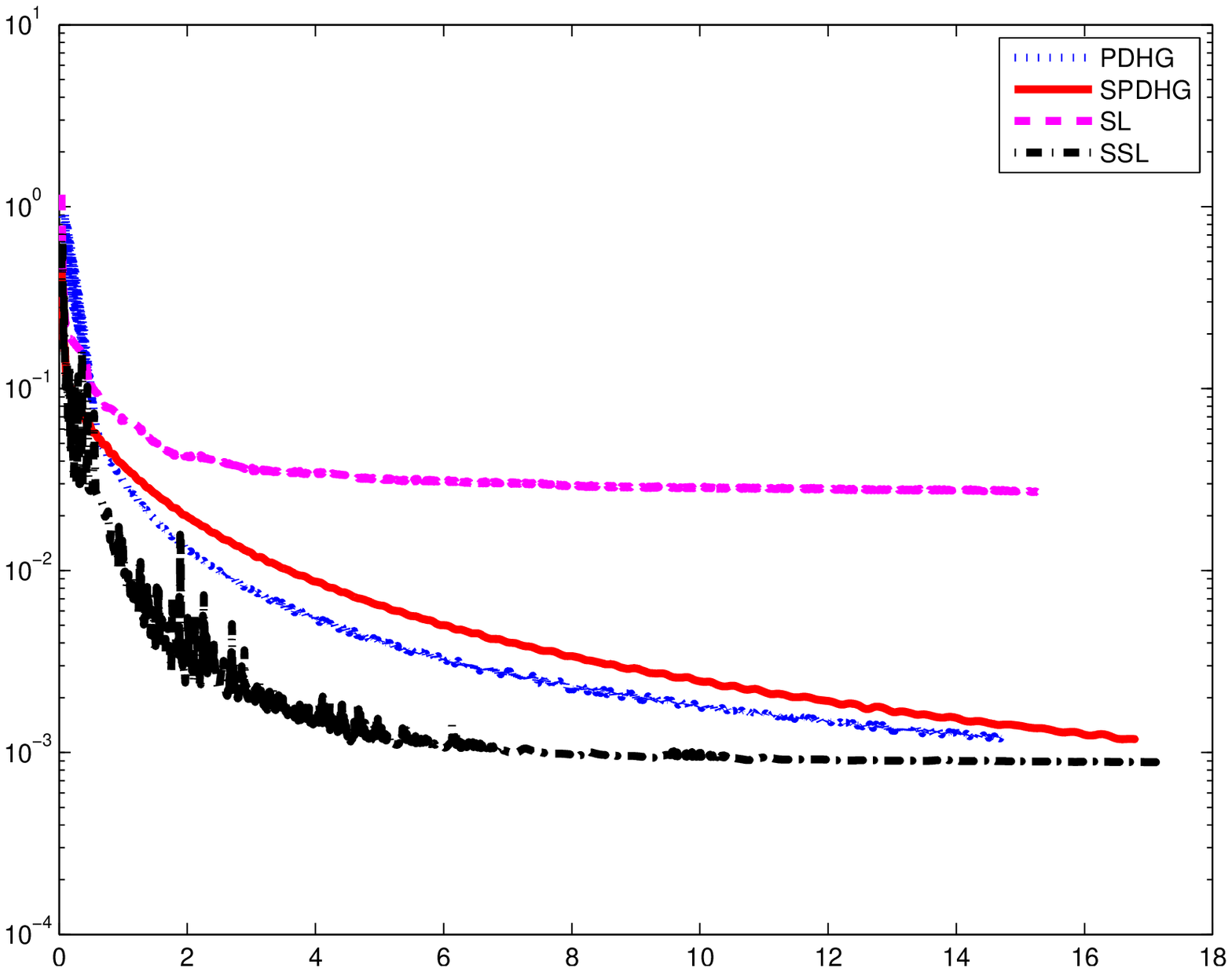}&\includegraphics[scale=\factor]{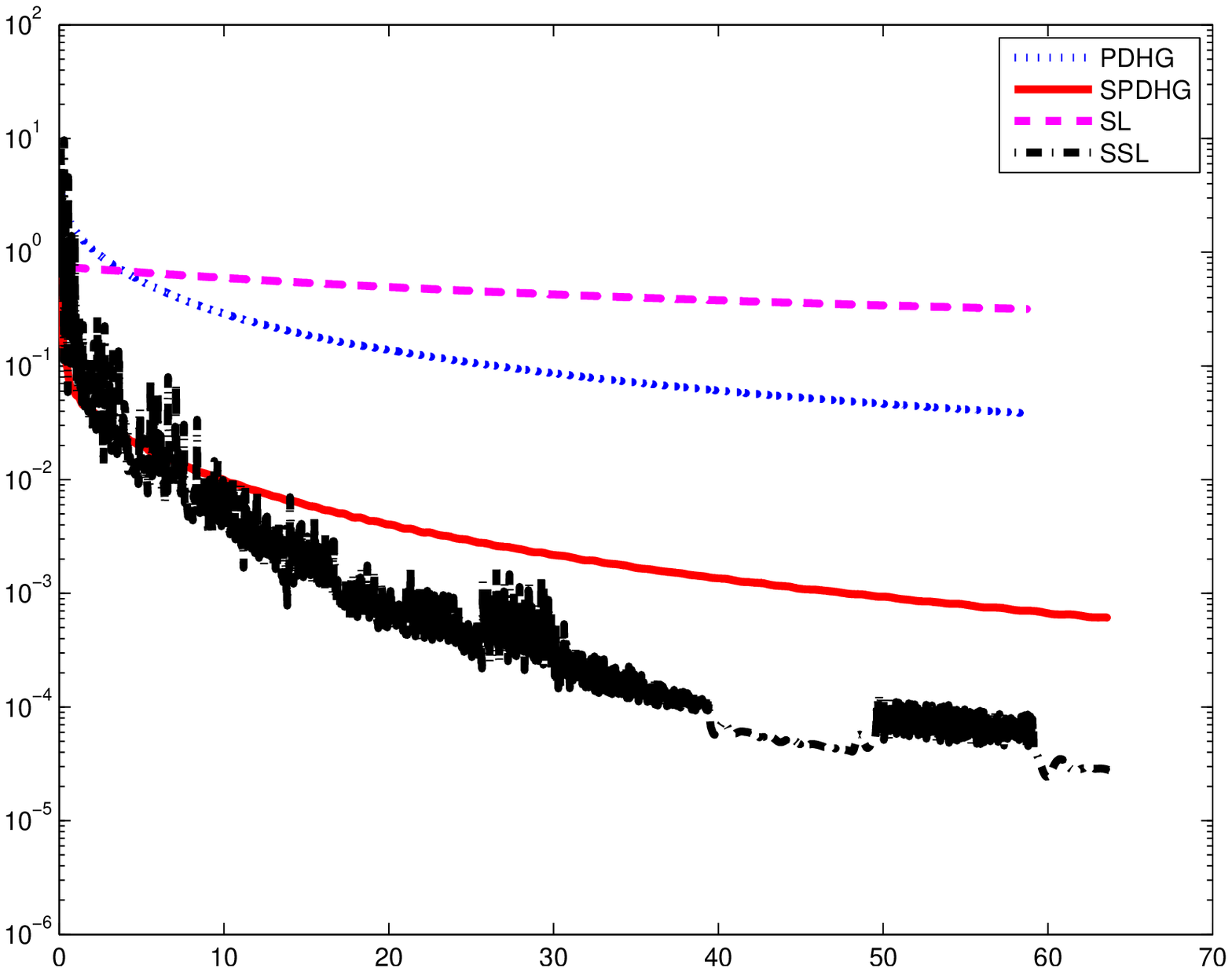}
\end{tabular}
\end{center}
\caption{Image deblurring results. Upper row: plots of the relative minimization error $e^k$. Lower row: plots of the relative difference from the optimal function value $f^k$. Left column: \emph{cameraman}. Middle column: \emph{micro}. Right column: \emph{phantom}. All plots are with respect to the computational time in seconds and use a logarithmic scale on the vertical axis.}\label{fig:1}
\end{figure}
The plots in Figure \ref{fig:1} have been obtained by running 3000 iterations of the algorithms, reporting the errors $e^k$, $f^k$ with respect to the computational time in seconds.

From the numerical experience we observe that the presence of the scaling can help to accelerate the progress towards the solution, with both stepsize selection strategies. As concerns the scaling matrix bounds, the best results are obtained by selecting large initial values for $\gamma_k$ (see Table \ref{table:1}) and, thus, for $L_k$, allowing more freedom to choose the scaling matrix especially at the first iterations.

It is also interesting to observe that the adaptive computation of $\ak$ combined with the proposed scaling technique seems to work quite well, leading to performances that are, in some cases, close to the `best' ones obtainable by manually tuning the stepsize sequence. Indeed, the PDHG and SPDHG performances are sensitive to the choice of $\tau_k$, and especially of $\ak$, making difficult to devise a general rule to select these sequences.

Comparing the first and the second row in Figure \ref{fig:1}, we also observe that a faster approaching to the solution $\x^*$ does not always correspond to a faster decrease of the objective function: we suppose that this phenomenon is due to the ill conditioning of problem \eqref{minsum}. For completeness, we experimentally observed that too large initial values of the primal stepsize $\ak$ in PDHG and SPDHG may produce an unbounded sequence $\{\uk\}$ and, as a consequence, the algorithms fail to converges: this indicates that the assumption on the $\epsilon$-subgradient boundedness is crucial.
\section{Conclusions}\label{sec7}
In this paper we proposed a generalization of the $\epsilon$-subgradient projection method with variable scaling matrix for nonsmooth, convex, constrained optimization, developing the related convergence analysis when the stepsize parameters are either provided as {\it a priori} selected sequences or dynamically computed by an adaptive procedure. Exploiting the duality principle, we described a special case of the proposed method which applies to the minimization of the sum of two convex functions. For a specific problem of this form in the image restoration framework, we fully detailed the algorithm, also suggesting a strategy to compute the scaling matrix. The numerical experience shows that the presence of a suitable variable scaling matrix can accelerate the progress of the iterates towards the solution. Moreover, the results obtained combining the variable scaling with the adaptive procedure for the computation of the stepsize parameter are encouraging.

Future work will be addressed to further investigate dynamic choices of the stepsize and of the scaling matrix, with the aim to devise effective `black--box' algorithms which are able to handle practical applications with a minimum of user supplied parameters.
\bibliography{scaling_subgradient}
\bibliographystyle{siam}
\end{document}